\documentclass[reqno]{amsart}
\usepackage{amssymb}
\usepackage{mathrsfs}
\usepackage{amsmath,amstext,amsfonts,verbatim}

\usepackage[latin1]{inputenc}

\usepackage{amsmath,amsthm,amssymb}

\usepackage{amsfonts}

\usepackage{amsmath,amssymb}
\usepackage{graphicx}
\usepackage{color}
\usepackage{indentfirst}

\theoremstyle{definition}

\theoremstyle{remark}

\numberwithin{equation}{section}

\newcommand{\neweq}[1]{\begin{equation}\label{#1}}

\def\phi{\varphi}

\def\incep{\left\{\begin{array}{cl} }
 \def\termin{\end{array}\right. }
\def\2af{2^*_\alpha}

\begin{document}

\title {A note on weak Banach mean equicoontinuity}{\textbf{******}}

\author{Zhongxuan Yang$^*$}
\address{College of Mathematics and Statistics, Chongqing University, Chongqing 401331, China}
\thanks{$^*$ Corresponding author.}
\email{20220601008@stu.cqu.edu.cn}
\thanks{The research was supported by NSF of China (No. 11671057) and NSF of Chongqing (Grant No. cstc2020jcyj-msxmX0694).}

\author{Xiaojun Huang}
\address{College of Mathematics and Statistics, Chongqing University, Chongqing 401331, China}
\email{hxj@cqu.edu.cn}

\keywords{weak mean equicontinuity, weak Banach mean equicontinuity, uniformly generic point, unique ergodic, uniform time average}

\subjclass[2010]{37A35, 37B40}

\date{}

\begin{abstract}
Consider a topological dynamical system $(X, T)$ endowed with the metric $d$. We introduce a novel function as $\overline{BF}(x, y) = \limsup_{n-m \rightarrow +\infty} \inf_{\sigma \in S_{n,m}} \frac{1}{n-m} \sum_{k=m}^{n-1} d\left(T^{k} x, T^{\sigma(k)} y\right)$, where the permutation group $S_{n,m}$ is utilized. It is demonstrated that $BF(x, y)$ exists when $x, y \in X$ are uniformly generic points.
Leveraging this function, we introduce the concept of weak Banach mean equicontinuity and establish that the dynamical system $(X, T)$ exhibits weak Banach mean equicontinuity if and only if the uniform time averages $f_B^{*}(x) = \lim_{n-m \rightarrow +\infty} \frac{1}{n-m} \sum_{k=m}^{n-1} f\left(T^{k} x\right)$ are continuous for all $f \in C(X)$. Finally, we demonstrate that in the case of a transitive system, the equivalence between weak Banach mean equicontinuity and weak mean equicontinuity is established.

\end{abstract}

\maketitle

\section{Introduction}

Let $X $ be a compact metric space with a metric $d$, and $T$ is a continuous map from $X$ to itself, we call the  pair $(X,T)$ a topological dynamical system.

It is a well-established fact that equicontinuous systems exhibit straightforward dynamical behaviors, with a primary focus on orbit analysis in the study of dynamical systems. When investigating long-term behaviors, particular emphasis has traditionally been placed on equicontinuous systems due to their inherent simplicity. However, the statistical properties of long-term behaviors are intricately tied to the cumulative effects of points within orbits. Consequently, it is justifiable to disregard the specific positions of certain points in orbits when delving into the statistical aspects of long-term behaviors. To address this consideration, the concept of mean-$L$-stable systems has been introduced in the literature [3,5,6]. A dynamical system denoted as $(X, T)$ is classified as mean-$L$-stable (refer to [5]) if, for any $\varepsilon > 0$, there exists a $\delta > 0$ such that $d(x, y) < \delta$ implies $d\left(T^{n} x, T^{n} y\right) < \varepsilon$ for all $n \in \mathbb{N}$, excluding a set of upper density less than $\varepsilon$. Notably, in recent contributions by Li, Tu, and Ye [1], the novel concepts of mean equicontinuous and Banach mean equicontinuous systems have been introduced. A dynamic system is termed mean equicontinuous (or Banach mean equicontinuous) if, for any $\varepsilon > 0$, there exists a $\delta > 0$ such that whenever $x, y \in X$ with $d(x, y) < \delta$,
$$
\limsup _{n \rightarrow+\infty} \frac{1}{n} \sum_{k=0}^{n-1} d\left(T^{k} x, T^{k} y\right)<\varepsilon .
$$
or
$$
\limsup _{n-m \rightarrow+\infty} \frac{1}{n-m} \sum_{k=m}^{n-1} d\left(T^{k} x, T^{k} y\right)<\varepsilon .
$$

In [1], the authors established a compelling result, demonstrating that a transitive dynamical system exhibits a dichotomy, manifesting as either almost mean equicontinuity or mean sensitivity. Furthermore, they rigorously established that a dynamical system achieves mean equicontinuity if and only if it attains the status of being mean-$L$-stable. The phenomenon of mean equicontinuity, a subject of heightened scholarly interest in recent years, has emerged as a pivotal area of study. Its significance is underscored by its profound connections with the ergodic properties inherent in measurable dynamical systems. For a comprehensive exploration of the intricacies surrounding mean equicontinuity and its correlated domains, interested scholars are encouraged to delve into the detailed survey presented in [1]. This survey provides a nuanced elucidation of mean equicontinuity, accompanied by further characterizations and insights elucidated in a multitude of references [7-17].

The classification of dynamical systems under the notion of mean equicontinuity is established through the utilization of the Besicovitch pseudometric, a topological construct stemming from the intrinsic $\overline{d}$ pseudometric prevalent in dynamic systems. In the seminal work [20], García-Ramos and Kwietniak undertook a pivotal modification by substituting the Besicovitch pseudometric with the Feldman-Katok pseudometric within the frameworks of mean equicontinuity and mean sensitivity. This adaptation gave rise to the introduction of Feldman-Katok continuity and a novel concept termed FK-sensitivity. The authors demonstrated that, for a minimal dynamical system, it inherently exhibits either Feldman-Katok continuity or FK-sensitivity. It is noteworthy that Feldman-Katok continuity, as juxtaposed with mean equicontinuity, manifests itself with a discernible level of weakened constraints. Moreover, in the scrutiny of the Feldman-Katok pseudometric, peculiar allowances for 'time delay' and 'space jump' are afforded. This peculiar feature emerges from the deliberate omission of point synchronization within orbit segments' distance, where only order preservation is mandated. Such an approach proves judicious when delving into the statistical intricacies of long-term dynamical behaviors.

In a recent study by Zheng and Zheng [14], the intricacies associated with the temporal order of points within orbits were deliberately overlooked. This strategic omission is grounded in the recognition that such order considerations lack meaningful significance in the exploration of statistical properties over extended temporal durations. In lieu of prioritizing point order, the authors introduced a novel mathematical construct known as the weak-mean pseudometric $\overline{F}$. This pseudometric, denoted by $\overline{F}$, was systematically formulated through the utilization of the permutation group $S_{n}$. For any  $x, y \in X$, we define
$$
\begin{aligned}
	\overline{F}(x, y) & =\limsup _{n \rightarrow+\infty} \inf _{\sigma \in S_{n}} \frac{1}{n} \sum_{k=1}^{n} d\left(T^{k} x, T^{\sigma(k)} y\right), \\
	\underline{F}(x, y) & =\liminf _{n \rightarrow+\infty} \inf _{\sigma \in S_{n}} \frac{1}{n} \sum_{k=1}^{n} d\left(T^{k} x, T^{\sigma(k)} y\right), \\
	N(\overline{F}) & =\{(x, y) \in X \times X \mid \overline{F}(x, y)=0\} \\
	N(\underline{F}) & =\{(x, y) \in X \times X \mid \underline{F}(x, y)=0\}
\end{aligned}
$$
where  $S_{n}$  is the n-order permutation group.  

{\bf Definition 1.1}(see [14])
Let  $(X, T)$  be a topological dynamical system. We say  $(X, T)$  is  $\overline{F}$-continuous at  $x \in X$  if for any  $\varepsilon>0$, there is a  $\delta>0$  such that whenever $ d(x, y)<\delta$, we have  $\overline{F}(x, y)<\varepsilon$. Denote by  $C(\overline{F})$  all  $\overline{F}$  continuous points. If  $C(\overline{F})=X$, we say  $(X, T)$  is  $\overline{F}$-continuous. In this case, we also call  $(X, T)$  weak mean equicontinuous.

{\bf Definition 1.2} (see [14])
Let  $(X, T)$  be a topological dynamical system.
 We say  $(X, T)$  is  $F$-continuous at  $x \in X$  if for any  $\varepsilon>0 $, there is a  $\delta>0$  such that whenever  $d(x, y)<\delta$, $ F(x, y)$  exists and $ F(x, y)<\varepsilon $. Denote by  $C(F)$  all  $F$-continuous points. If  $C(F)=X $, we say  $(X, T)$  is  $F$-continuous. 

In their seminal work [14], the authors established that weak mean equicontinuous systems precisely coincide with those wherein the time-average operator successfully maintains the continuity of observed functions. It is noteworthy that the weak-mean pseudometric $\overline{F}$ is equivalent to pseudometric $\widetilde{d}_{FK}$ as defined in [19]. The pseudometric $\widetilde{d}_{FK}$ is characterized by the omission of the order in the definition of the Feldman-Katok pseudometric. Consequently, it follows that  weak mean equicontinuity is weaker than  Feldman-Katok continuity.

In the seminal work [1], the concept of Banach mean equicontinuity was introduced, prompting an inquiry into the equivalence between Banach mean equicontinuity and mean equicontinuity in minimal dynamical systems. This fundamental question found an affirmative resolution in subsequent research, precisely in [18]. Further contributing to this discourse, the authors in [7] extended their investigation, demonstrating the broader result that the equivalence of Banach mean equicontinuity and mean equicontinuity holds universally for arbitrary dynamical systems.

In pursuit of articulating a lucid exposition of our findings, we establish the following notational framework. For arbitrary elements $x, y \in X$ within the dynamical system under consideration, we proceed to define
$$
\begin{aligned}
	\overline{BF}(x, y) & =\limsup _{n-m \rightarrow+\infty} \inf _{\sigma \in S_{n,m}} \frac{1}{n-m} \sum_{k=m}^{n-1} d\left(T^{k} x, T^{\sigma(k)} y\right), \\
	\underline{BF}(x, y) & =\liminf _{n-m \rightarrow+\infty} \inf _{\sigma \in S_{n,m}} \frac{1}{n-m} \sum_{k=m}^{n-1} d\left(T^{k} x, T^{\sigma(k)} y\right), \\
	N(\overline{BF}) & =\{(x, y) \in X \times X \mid \overline{BF}(x, y)=0\} \\
	N(\underline{BF}) & =\{(x, y) \in X \times X \mid \underline{BF}(x, y)=0\}
\end{aligned}
$$
where  $S_{n,m}$  is the  permutation group of $\{m,m+1,\cdots,n-1\}$.

If  $\overline{BF}(x, y)=\underline{BF}(x, y)$, we say  $BF(x, y)$  exists, and define
$$
\begin{aligned}
	BF(x, y)=\lim _{n-m \rightarrow+\infty} \inf _{\sigma \in S_{n,m}} \frac{1}{n-m} \sum_{m}^{n-1} d\left(T^{k} x, T^{\sigma(k)} y\right), \\
	N(BF)=\{(x, y) \in X \times X \mid BF(x, y)=0\} .
\end{aligned}
$$
It is easy to obtain that
$
N(\overline{BF})=N(BF) \subset N(\underline{BF}).$

{\bf Definition 1.3}
Let  $(X, T)$  be a topological dynamical system. We say  $(X, T)$  is  $\overline{BF}$-continuous at  $x \in X$  if for any  $\varepsilon>0$, there is a  $\delta>0$  such that whenever $ d(x, y)<\delta$, we have  $\overline{BF}(x, y)<\varepsilon$. Denote by  $C(\overline{BF})$  all  $\overline{BF}$-continuous points. If  $C(\overline{BF})=X$, we say  $(X, T)$  is  $\overline{BF}$-continuous. In this case, we also call  $(X, T)$  weak  Banach mean equicontinuous.

{\bf Definition 1.4}  Let  $(X, T)$  be a topological dynamical system. We say  $(X, T)$  is  $BF$-continuous at  $x \in X$  if for any  $\varepsilon>0 $, there is a  $\delta>0$  such that whenever  $d(x, y)<\delta$, $ BF(x, y)$  exists and $ BF(x, y)<\varepsilon $. Denote by  $C(BF)$  all  $BF$-continuous points. If  $C(BF)=X $, we say  $(X, T)$  is  $BF$-continuous.

The functions $\overline{F}(x, y)$, $\overline{F}(x, y)$, $\overline{BF}(x, y)$ and $\underline{BF}(x, y)$ serve as pseudometrics, quantifying the distinctions between the distributions of $\operatorname{Orb}(x,T)$ and $\operatorname{Orb}(y,T)$.
Here, $\operatorname{Orb}(x,T)=\left\{x, T x, T^{2} x, \cdots\right\}$ and $\operatorname{Orb}(y,T)=\left\{y, T y, T^{2} y, \cdots\right\}$ denote the orbits of $x$ and $y$.

In [14], Zheng and Zheng demonstrated that when $x$ and $y$ are generic points (refer to Section 2), the equality $\overline{F}(x, y)=\underline{F}(x, y)$ holds. Furthermore, in [18], Downarowicz and Glasner introduced the concept of uniformly generic points (see Section 2), which prompted a natural inquiry: if $x$ and $y$ are uniformly generic points, does the equality $\overline{BF}(x, y)=\underline{BF}(x, y)$ persist?

Moving to the realm of ergodicity, Fomin [6] establishes that a minimal mean-$L$-stable system is uniquely ergodic. Building upon this, Zheng and Zheng [14] extend the discourse to transitive weak mean equicontinuous systems, proving their unique ergodicity. The equivalences established in their work are noteworthy: (1) $(X, T)$ is uniquely ergodic; (2) $N(F)=X \times X$; (3) $N(\underline{F})=X \times X$. A compelling question emerges: If we substitute weak mean equicontinuity with weak Banach mean equicontinuity in the context of dynamical systems, what implications would unfold?

Indeed, it is evident that Banach mean equicontinuity implies weak Banach mean equicontinuity, and similarly, weak Banach mean equicontinuity implies weak mean equicontinuity. However, it is noteworthy that in the general case, weak Banach mean equicontinuity does not necessarily imply Banach mean equicontinuity. 

On one hand,  $BF$-continuity implies  $\overline{BF}$-continuity. Conversely, within an $\overline{BF}$-continuous system denoted as $(X, T)$, we can rigorously demonstrate that all points therein are uniformly generic points.  This significant result, coupled with the insights from Theorem 4.2 and Other types of equicontinuities. In summary, the relationships can be succinctly summarized as follows:

equicontinuity  $\Rightarrow$  mean equicontinuity  $\Leftrightarrow$  Banach mean equicontinuity  $\Rightarrow$  weak Banach mean equicontinuity($\overline{BF}$-continuity) $\Leftrightarrow$ $BF$-continuity.

Banach mean equicontinuity  $\Leftrightarrow$   mean equicontinuity  $\Rightarrow$  weak  mean equicontinuity($\overline{F}$-continuity)$\Leftrightarrow$ $F$-continuity.

We designate a dynamical system as chaotic when the trajectories of points in its orbits exhibit sensitivity to initial conditions. Despite the potential chaotic behavior observed in weak Banach mean equicontinuous systems, a nuanced perspective grounded in measure theory reveals their inherent stability. This stability emanates from the fact that the probability distributions characterizing the positions of points within orbits demonstrate insensitivity to variations in initial values.

The Birkhoff Ergodic Theorem establishes the integrability of time averages $f^{*}$(refer to Section 2) for any integrable function $f$. The inquiry naturally arises: under what conditions can the time average operator ensure the continuity of observed functions? Zheng and Zheng [14] provide a comprehensive answer, revealing that $(X, T)$ is weak mean equicontinuous if and only if the time averages $f^{*}$ remain continuous for all $f \in C(X)$. Motivated by the concept of uniformly generic points, we introduce the notion of uniformly time average $f^{*}_{B}$ (refer to Section 2). Our focus shifts to discerning whether $(X, T)$ being weak Banach mean equicontinuous implies the continuity of uniformly time averages $f^{*}_{B}$ for all $f \in C(X)$. Building upon the groundwork laid by Qiu and Zhao [20], who demonstrated the equivalence of Banach mean equicontinuity and mean equicontinuity, a natural extension emerges: the investigation into the conditions under which weak Banach mean equicontinuity aligns with weak mean equicontinuity.

This paper is structured as follows: In Section 2, we present fundamental concepts and results essential to our investigation. Section 3 establishes key propositions concerning $\overline{BF}$ and $N(BF)$, pivotal to subsequent discussions. Section 4 is dedicated to the rigorous proof of the existence of $BF(x, y)$ under the condition that $x$ and $y$ are uniformly generic points. Transitioning to Section 5, we unravel the unique ergodicity of a transitive weak Banach mean equicontinuous system. Finally, Section 6, we show that for all the points of a  weak Banach mean equicontinuous system are uniformly generic points. Specifically, a system is categorized as weak Banach mean equicontinuous if and only if uniformly time averages exhibit continuity. For transitive systems, the equivalence between weak Banach mean equicontinuity and weak mean equicontinuity is established.

\section{Preliminaries}

In this section we recall some notions and results of topological dynamical system.  which are needed in our paper. Note that  $\mathbb{N}$  denotes the set of all non-negative integers and  $\mathbb{N}^{+}$ denotes the set of all positive integers in this paper.

{\bf 2.1} Let  $F \subset \mathbb{N}$, we define the upper density  $\overline{D}(F)$  of  $F$  by
$$
\overline{D}(F)=\limsup _{n \rightarrow+\infty} \frac{\#(F \cap[0, n-1])}{n},
$$
where  $\#(\cdot)$  is the number of elements of a set.

Similarly, the lower density  $\underline{D}(F)$  of $ F$  is defined by
$$
\underline{D}(F)=\liminf _{n \rightarrow+\infty} \frac{\#(F \cap[0, n-1])}{n} .
$$
We call  $ F$  has density $ D(F)$  if  $\overline{D}(F)=\underline{D}(F) $.

The upper Banach density  $\overline {BD}(F)$  is defined by
$$
\overline {BD}(F)=\limsup _{n-m \rightarrow \infty} \frac{\#(F \cap[m, n-1])}{N-M} .$$

{\bf 2.2} Suppose  $(X, T)$  is a topological dynamical system. The  $\sigma$ -algebra of Borel subsets of  $X$  will be denoted According to  $\mathscr{B}(X) $. Let  $M(X) $ be the collection of all regular Borel probability measures defined on the measurable space  $(X, \mathscr{B}(X))$. In the weak* topology, $ M(X)$  is a nonempty compact set.

We say  $\mu \in M(X)$  is  $T$-invariant if  $\mu\left(T^{-1}(A)\right)=\mu(A)$  holds for any  $A \in \mathscr{B}(X) $. Denote According to $ M(X, T)$  the collection of all $ T $-invarant regular Borel probability measures defined on the measurable space  $(X, \mathscr{B}(X))$. In the weak  ${ }^{*} $ topology,  $M(X, T)$  is a nonempty compact convex set.

We say  $\mu \in M(X, T)$  is ergodic if for any  $A \in \mathscr{B}(X)$  with  $T^{-1} A=A$, $ \mu(A)=0  $ or $ \mu(A)=1 $ holds. Denote by  $E(X, T)$  the collection of all ergodic measures on $ (X, T)$. As well known, $ E(X, T)$  is the collection of all extreme points of  $M(X, T)$ and  $E(X, T)$  is nonempty. 

We say  $(X, T)$  is uniquely ergodic if  $E(X, T)$  is singleton. Since  $E(X, T)$  is the set of extreme points of  $M(X, T)$, then  $(X, T)$  is unique ergodic if and only if  $M(X, T) $ is singleton. For  $\mu \in M(X, T)$, the support of  $\mu$  is defined by  $\operatorname{supp}(\mu)=\{x \in X: \mu(U)>  0 \text{\ for any neighborhood\ }  U  \text{\ of\ }  x\} $. A measure  $\mu$  on  $X$  has full support if  $\operatorname{supp}(\mu)=X $. It is known that the support of an ergodic measure is a transitive subsystem. The support of a topological dynamical system  $(X, T)$  is defined by  $\operatorname{supp}(X, T)=\overline {\bigcup \{\operatorname{supp}(\mu): \mu \in M(X, T)\}} $.

Given  $x \in X $, we have  $\left\{\frac{1}{n-m} \sum_{k=m}^{n-1} \delta_{T^{k} x}\right\}_{m\in\mathbb{N},n>m} \subset M(X)$, where $ \delta_{x} $ is the Dirac measure supported on $x$. Denote by $ M_{B,x}$  the collection of all limit points of  $\left\{\frac{1}{n-m} \sum_{k=m}^{n-1} \delta_{T^{k} x}\right\}_{m\in\mathbb{N},n>m} .$  Since  $M(X) $ is compact, we have  $M_{B,x} \neq \emptyset$. Moreover, $ M_{B,x} \subset M(X, T)$. We call  $M_{B,x}$  the measure set generated by  $x$.

A point $ x \in X $ is called  generic point (see [18]) if for any  $f \in C(X)$, the time average
$$
f^{*}(x)=\lim _{n \rightarrow+\infty} \frac{1}{n} \sum_{k=1}^{n} f\left(T^{k} x\right)
$$
exist.

A point $ x \in X $ is called uniformly generic point (see [18]) if for any  $f \in C(X)$, the uniform time average
$$
f_B^{*}(x)=\lim _{n-m \rightarrow+\infty} \frac{1}{n-m} \sum_{k=m}^{n-1} f\left(T^{k} x\right)
$$
exists, let  $Q$  denote the set of all uniformly generic points. It is easy to derive that  $x$  is a uniformly generic point if and only if  $M_{B,x} $ consists of a single measure. We call  $\mu \in M_{B,x}$  is generated by  $x$  if  $x$  is a uniformly generic point.  We call a uniformly generic point  $x$  is an ergodic point if the invariant measure generated by  $x$  is ergodic.

A Borel subset  $E \subset X$  is said to have invariant measure one if  $\mu(E)=1$  for all $ \mu \in M(X, T) $. 

Next, we define uniformly physical measures in a general way.

{\bf Definition 2.1} Let  $(X, T)$  be a topological dynamical system and  $m \in M(X)$. We say  $\mu \in M(X, T)$  is a uniformly physical measure with respect to  m  if  $m(B(\mu))>0$, where
$$
B(\mu)=\{x \in X \mid \lim _{n-m \rightarrow+\infty} \frac{1}{n-m} \sum_{k=m}^{n-1} \delta_{T^{k} x}=\mu\}
$$
For any  $A \subset X $, let
$$
\chi_{A}(x)=\left\{\begin{array}{ll}
	1, & x \in A \\
	0, & x \notin A
\end{array}\right.
$$
The following Lemma is well known (see [2], page:149).

{\bf Lemma 2.2} 
 Let  $(X, T)$  be a topological dynamical system. If  $x \in X$  is a  uniformly generic point and  $\mu $ is generated by $x $, then for any open set  $U \subset X$  and any closed set  $V \subset X $, we have
$$
\liminf _{n-m \rightarrow+\infty} \frac{1}{n-m} \sum_{k=m}^{n-1} \chi_{U}\left(T^{k} x\right) \geq \mu(U) \text { and } \limsup _{n-m \rightarrow+\infty} \frac{1}{n-m} \sum_{k=m}^{n-1} \chi_{V}\left(T^{k} x\right) \leq \mu(V) .
$$

Given  $\mu \in M(X) $. Since  $X $ is compact, there are finite mutually disjoint subsets of  $X$  such that the diameter of each subset is small enough and the sum of their measures are closed enough to one. Hence we have the following result.

{\bf Lemma 2.3}  Let  $\mu \in M(X)$. Then for any  $\varepsilon>0$ , there are finite mutually disjoint closed sets 
$ \left\{\Lambda_{k}\right\}_{k=1}^{k_{0}} $ such that
$$
\mu\left(\bigcup_{k=1}^{k_{0}} \Lambda_{k}\right) \geq 1-\varepsilon \text { and } \operatorname{diam}\left(\Lambda_{k}\right) \leq \varepsilon, \forall k=1,2, \cdots, k_{0} .
$$
Similarly, there are finite mutually disjoint open sets  $\left\{V_{s}\right\}_{s=1}^{s_{0}} $ such that
$$
\mu\left(\bigcup_{s=1}^{s_{0}} V_{s}\right) \geq 1-\varepsilon \text { and } \operatorname{diam}\left(V_{s}\right) \leq \varepsilon, \forall s=1,2, \cdots, s_{0}
$$
\section{Some basic properties of  $BF(x, y)$  and  $N(BF)$ }
In this section, we will show some properties of  $\overline{BF}(x, y) $ and  $N(BF)$, which play an important role in the following sections.

{\bf Proposition 3.1} Let  $(X, T)$  be a topological dynamical system. Then

(1) For any sequences  $\left\{x_{k}\right\}_{k=m}^{n-1}$  and  $\left\{y_{k}\right\}_{k=m}^{n-1}$  of $ X$, we have
$$
\inf _{\sigma \in S_{n,m}} \sum_{k=m}^{n-1} d\left(x_{k}, y_{\sigma(k)}\right)=\inf _{\sigma \in S_{n,m}} \sum_{k=m}^{n-1} d\left(y_{k}, x_{\sigma(k)}\right)
$$
In particular, for any  $x, y \in X$, we have
$$
\inf _{\sigma \in S_{n,m}} \frac{1}{n-m} \sum_{k=m}^{n-1} d\left(T^{k} x, T^{\sigma(k)} y\right)=\inf _{\sigma \in S_{n,m}} \frac{1}{n-m} \sum_{k=m}^{n-1} d\left(T^{k} y, T^{\sigma(k)} x\right) .
$$

(2) For any sequences  $\left\{x_{k}\right\}_{k=m}^{n-1},\left\{y_{k}\right\}_{k=m}^{n-1}$  and $\left\{z_{k}\right\}_{k=m}^{n-1}$  of  $X$ , we have
$$
\inf _{\sigma \in S_{n,m}} \sum_{k=m}^{n-1} d\left(x_{k}, z_{\sigma(k)}\right) \leq \inf _{\sigma \in S_{n,m}} \sum_{k=m}^{n-1} d\left(x_{k}, y_{\sigma(k)}\right)+\inf _{\sigma \in S_{n,m}} \sum_{k=m}^{n-1} d\left(y_{k}, z_{\sigma(k)}\right).$$
In particular, for any  $x, y, z \in X$ , we have
$$
\inf _{\sigma \in S_{n,m}} \frac{1}{n-m} \sum_{k=m}^{n-1} d\left(T^{k} x, T^{\sigma(k)} z\right) \leq \inf _{\sigma \in S_{n,m}} \frac{1}{n-m} \sum_{k=m}^{n-1} d\left(T^{k} x, T^{\sigma(k)} y\right)+\inf _{\sigma \in S_{n,m}} \frac{1}{n-m} \sum_{k=m}^{n-1} d\left(T^{k} y, T^{\sigma(k)} z\right).
$$

(3) For any  $x, y \in X$ , we have
$$
\overline{BF}(x, y)=\overline{BF}(y, x)
$$
and
$$
\underline{BF}(x, y)=\underline{BF}(y, x).
$$

(4) For any  $x, y, z \in X$ , we have
$$
\overline{BF}(x, z) \leq \overline{BF}(x, y)+\overline{BF}(y, z)
$$
and
$$
\underline{BF}(x, z) \leq \underline{BF}(x, y)+\overline{BF}(y, z).
$$
{\bf Proof} (1) There exists a  $\sigma_{1} \in S_{n,m}$  such that
$$
\sum_{k=m}^{n-1} d\left(x_{k}, y_{\sigma_{1}(k)}\right)=\inf _{\sigma \in S_{n,m}} \sum_{k=m}^{n-1} d\left(x_{k}, y_{\sigma(k)}\right) \text {. }
$$
Let  $\sigma_{2} \in S_{n,m}$  such that  $\sigma_{1} \sigma_{2}=\sigma_{2} \sigma_{1}$  be the indentity element of  $S_{n,m}$ . Then we have
$$
\sum_{k=m}^{n-1}d\left(x_{k}, y_{\sigma_{1}(k)}\right)=\sum_{k=m}^{n-1} d\left(x_{\sigma_{2} \sigma_{1}(k)}, y_{\sigma_{1}(k)}\right)=\sum_{k=m}^{n-1} d\left(x_{\sigma_{2}(k)}, y_{k}\right) \geq \inf _{\sigma \in S_{n,m}} \sum_{k=m}^{n-1} d\left(y_{k}, x_{\sigma(k)}\right) .
$$
Thus,
$$
\inf _{\sigma \in S_{n,m}} \sum_{k=m}^{n-1}d\left(x_{k}, y_{\sigma(k)}\right) \geq \inf _{\sigma \in S_{n,m}} \sum_{k=m}^{n-1} d\left(y_{k}, x_{\sigma(k)}\right)
$$
Similary, we can get
$$
\inf _{\sigma \in S_{n,m}} \sum_{k=m}^{n-1} d\left(y_{k}, x_{\sigma(k)}\right) \geq \inf _{\sigma \in S_{n,m}} \sum_{k=m}^{n-1} d\left(x_{k}, y_{\sigma(k)}\right)
$$
Hence,
$$
\inf _{\sigma \in S_{n,m}} \sum_{k=m}^{n-1} d\left(x_{k}, y_{\sigma(k)}\right)=\inf _{\sigma \in S_{n,m}} \sum_{k=m}^{n-1} d\left(y_{k}, x_{\sigma(k)}\right)
$$
(2) There are  $\sigma_{1}, \sigma_{2} \in S_{n,m}$  such that
$$
\sum_{k=m}^{n-1}d\left(x_{k}, y_{\sigma_{1}(k)}\right)=\inf _{\sigma \in S_{n,m}} \sum_{k=m}^{n-1} d\left(x_{k}, y_{\sigma(k)}\right)
$$
and
$$
\sum_{k=m}^{n-1} d\left(y_{k}, z_{\sigma_{2}(k)}\right)=\inf _{\sigma \in S_{\mathrm{n,m}}} \sum_{k=m}^{n-1} d\left(y_{k}, z_{\sigma(k)}\right).
$$
Let  $\sigma_{3}=\sigma_{2} \sigma_{1}$, then we have
$$
\begin{aligned}
	\sum_{k=m}^{n-1} d\left(x_{k}, z_{\sigma_{3}(k)}\right) & \leq \sum_{k=m}^{n-1} d\left(x_{k}, y_{\sigma_{1}(k)}\right)+\sum_{k=m}^{n-1} d\left(y_{\sigma_{1}(k)}, z_{\sigma_{3}(k)}\right) \\
	& =\sum_{k=m}^{n-1} d\left(x_{k}, y_{\sigma_{1}(k)}\right)+\sum_{k=m}^{n-1} d\left(y_{k}, z_{\sigma_{2}(k)}\right) \\
	& =\inf _{\sigma \in S_{n,m}} \sum_{k=m}^{n-1} d\left(x_{k}, y_{\sigma(k)}\right)+\inf _{\sigma \in S_{n,m}} \sum_{k=m}^{n-1} d\left(y_{k}, z_{\sigma(k)}\right)
\end{aligned}
$$
Thus,
$$
\inf _{\sigma \in S_{n,m}} \sum_{k=m}^{n-1} d\left(x_{k}, z_{\sigma(k)}\right) \leq \inf _{\sigma \in S_{n,m}} \sum_{k=m}^{n-1} d\left(x_{k}, y_{\sigma(k)}\right)+\inf _{\sigma \in S_{n,m}} \sum_{k=m}^{n-1} d\left(y_{k}, z_{\sigma(k)}\right) .
$$
According to (1) and (2), we can easily deduce (3) and (4).

$\overline{BF}(x, y)$  and  $\underline{BF}(x, y)$  are functions which can measure the difference between the distributions of  $\operatorname{Orb}(x)$  and  $\operatorname{Orb}(y)$. And if $ x$  and  $y$  are in the same orbit, the distributions of  $\operatorname{Orb}(x)$  and $ \operatorname{Orb}(y)$  are same. Therefore,  $\overline{BF}(x, y)=BF(x, y)=0$. 

Next, we will indicate that $ \left.\overline{BF}\right|_{\operatorname{Orb}(x) \times \operatorname{Orb}(y)} $ and $ \underline{BF} |_{\operatorname{Orb}(x) \times \operatorname{Orb}(y)}$  are constant for any  $x, y \in X .$

{\bf Proposition 3.2} Let  $(X, T)$  be a topological dynamical system. For any  $x, y \in X$  and  $r, s \in \mathbb{N}$ , we have
$$
\overline{BF}(T^{r} x, T^{s} y)=\overline{BF}(x, y)
$$
and
$$
\underline{BF}\left(T^{r} x, T^{s} y\right)=\underline{BF}(x, y)
$$
If $ BF(x, y)$  exists, we also have
$$
BF\left(T^{r} x, T^{s} y\right)=BF(x, y)
$$
{\bf Proof} According to Proposition 3.1, we have
$$
\underline{BF}\left(T^{r} x, y\right) \leq \underline{BF}(x, y)+\overline{BF}\left(T^{r} x, x\right)=\underline{BF}(x, y) .
$$
On the other hand, we have
$$
\underline{BF}(x, y) \leq \underline{BF}\left(T^{r} x, y\right)+\overline{BF}\left(x, T^{r} x\right)=\underline{BF}\left(T^{r} x, y\right) .
$$
Thus,
$$
\underline{BF}\left(T^{r} x, y\right)=\underline{BF}(x, y) .
$$
Similarly, we can deduce that
$$
\underline{BF}\left(T^{r} x, T^{s} y\right)=\underline{BF}\left(T^{r} x, y\right).
$$
Hence, we have
$$
\underline{BF}\left(T^{r} x, T^{s} y\right)=\underline{BF}(x, y).
$$
Similarly, we have $$
\overline{BF}\left(T^{r} x, y\right) \leq \overline{BF}(x, y)+\overline{BF}\left(T^{r} x, x\right)=\overline{BF}(x, y) .
$$
On the other hand, we have
$$
\overline{BF}(x, y) \leq \overline{BF}\left(T^{r} x, y\right)+\overline{BF}\left(x, T^{r} x\right)=\overline{BF}\left(T^{r} x, y\right) .
$$
Thus,
$$
\overline{BF}\left(T^{r} x, y\right)=\overline{BF}(x, y) .
$$
Similarly, we can deduce that
$$
\overline{BF}\left(T^{r} x, T^{s} y\right)=\overline{BF}\left(T^{r} x, y\right).
$$
Hence, we have
$$
\overline{BF}\left(T^{r} x, T^{s} y\right)=\overline{BF}(x, y).
$$
Combining the above results, we can deduce that if $ BF(x, y)$  exists, then
$
BF\left(T^{r} x, T^{s} y\right)=BF(x, y).
$

According to Proposition 3.2, we can deduce  $N(BF)$  and  $N(\underline{BF})$  are both invariant sets with respect to  $T^{r} \times T^{s}$  for any $ r, s \in \mathbb{N}$. Given  $x \in X$ , let
$$
N(BF, x)=\{y \in X \mid BF(x, y)=0\}
$$
and
$$
N(\underline{BF}, x)=\{y \in X \mid \underline{BF}(x, y)=0\} .
$$
Then according to Proposition 3.2  we derive that  $N(BF, x)$  and  $N(\underline{BF}, x)$  are both invariant sets.

The following proposition provides a way to estimate the upper bound of  $\overline{BF}(x, y)$.

{\bf Proposition 3.3} Let  $(X, T)$  be a topological dynamical system, and  $\left\{U_{s}\right\}_{s=1}^{s_{0}}$  be mutually disjoint subsets of  $X$. Given $ x, y \in X $. If the following two conditions hold:

(1) There is  $\varepsilon>0$  such that $diam  \left(U_{s}\right) \leq \varepsilon$  holds for any  $s \in\left\{1,2, \cdots, s_{0}\right\} ;$

(2) For any  $s \in\left\{1,2, \cdots, s_{0}\right\}$ , there is  $a_{s} \geq 0$  such that
$$
\liminf _{n-m \rightarrow+\infty} \frac{1}{n-m}\sum_{k=m}^{n-1} \chi_{U_{s}}\left(T^{k} x\right) \geq a_{s}
$$
and
$$
\liminf _{n-m \rightarrow+\infty} \frac{1}{n-m}\sum_{k=m}^{n-1} \chi_{U_{s}}\left(T^{k} y\right) \geq a_{s}
$$
Then, we have
$$
\overline{BF}(x, y) \leq \varepsilon \sum_{s=1}^{s_{0}} a_{s}+M\left(1-\sum_{s=1}^{s_{0}} a_{s}\right)
$$
where  $M=\operatorname{diam}(X) .$

{{\bf Proof} Given  $\delta>0 $, there is an  $n_{1}>0 $ such that for any  $n-m>n_{1}$  and any  $s \in\left\{1,2, \cdots, s_{0}\right\} $, we have
	$$
	\frac{1}{n-m}\sum_{k=m}^{n-1} \chi_{U_{s}}\left(T^{k} x\right) \geq \liminf _{n \rightarrow+\infty} \frac{1}{n-m}\sum_{k=m}^{n-1} \chi_{U_{s}}\left(T^{k} x\right)-\delta \geq a_{s}-\delta .\ \ \ \ \ \ \text {(3.1)}
	$$
	Similarly, there is an  $n_{2}>0$  such that for any  $n-m>n_{2}$  and any  $s \in\left\{1,2, \cdots, s_{0}\right\} ,$ we have
	$$
	\frac{1}{n-m}\sum_{k=m}^{n-1} \chi_{U_{8}}\left(T^{k} y\right) \geq a_{s}-\delta       \ \ \ \ \ \ \ \ \ \ \ \  \ \  \ \text {(3.2)}
	$$
	Set  $n-m>\max \left\{n_{1}, n_{2}\right\}$. For any $ s \in\left\{1,2, \cdots, s_{0}\right\} $, let
	$$
	N_{n,m}\left(U_{s}, x\right)=\left\{k \in \mathbb{N}^{+} \mid T^{k} x \in U_{s}, k\in [m,n-1] \right\} \text { and } N_{n,m}\left(U_{s}, y\right)=\left\{k \in \mathbb{N}^{+} \mid T^{k} y \in U_{s}, k\in [m,n-1]\right\} .$$
	According to (3.1) and (3.2), there exists  $r_{s, n,m} \in \mathbb{N} $ with
	$$
	a_{s} (n-m)+1 \geq r_{s, n,m} \geq\left(a_{s}-\delta\right) (n-m)      \ \ \ \ \ \ \ \ \ \ \ \  \ \  \ \text {(3.3)}
	$$
	such that
	$$
	\#\left(N_{n.m}\left(U_{s}, x\right)\right) \geq r_{s, n,m} \text { and }  \#\left(N_{n,m}\left(U_{s}, y\right)\right) \geq r_{s, n,m} .
	$$
	Thus there are subsets $ \left\{k_{x, s, r}\right\}_{r=1}^{r_{s, n,m}}$  and $ \left\{k_{y, s, r}\right\}_{r=1}^{r_{s, n,m}}$  of $ \mathbb{N}^{+}$ such that
	$$
	\left\{k_{x, s, 1}, k_{x, s, 2}, \cdots, k_{x, s, r_{s, n,m}}\right\} \subset N_{n,m}\left(U_{s}, x\right)
	$$
	and
	$$
	\left\{k_{y, s, 1}, k_{y, s, 2}, \cdots, k_{y, s, r_{s, n,m}}\right\} \subset N_{n,m}\left(U_{s}, y\right) .
	$$
	Since  $N_{n,m}\left(U_{s_{1}}, x\right) \cap N_{n,m}\left(U_{s_{2}}, x\right)=\emptyset$  and  $N_{n,m}\left(U_{s_{1}}, y\right) \cap N_{n,m}\left(U_{s_{2}}, y\right)=\emptyset$  hold for any  $s_{1}, s_{2} \in\left\{1,2, \cdots, s_{0}\right\}$  with  $s_{1} \neq s_{2}$ , there is $ \sigma_{n,m} \in S_{n,m}$  such that
	$$
	\sigma_{n,m}\left(k_{x, s, r}\right)=k_{y, s, r}        \ \ \ \ \ \ \ \ \ \ \ \  \ \  \ \text {(3.4)}
	$$
	holds for any $ s \in\left\{1,2, \cdots, s_{0}\right\}$  and any $ r \in\left\{1,2, \cdots, r_{s, n,m}\right\} .$
	Let
	$$
	A=\bigcup_{s=1}^{s_{0}}\left\{k_{x, s, r}\right\}_{r=1}^{r_{s, n,m}} \text { and } B=\{m,m+1, \cdots, n-1\} \backslash A \text {. }
	$$
	According to (3.3),  we have
	$$
	\#(A)=\sum_{s=1}^{s_{0}} r_{s, n,m} \leq \sum_{s=1}^{s_{0}}\left(a_{s} (n-m)+1\right)     \ \ \ \ \ \ \ \ \ \ \ \  \ \  \ \text {(3.5)}
	$$
	and
	$$
	\#(B)=(n-m)-\sum_{s=1}^{s_{0}} r_{s, n,m} \leq n-m-\sum_{s=1}^{s_{0}}\left(a_{s}-\delta\right) (n-m) . \ \ \ \ \ \ \ \ \ \ \ \  \ \  \ \text {(3.6)}
	$$
	According to (3.4), we deduce that
	$$
	d\left(T^{k} x, T^{\sigma_{n,m}(k)} y\right) \leq \varepsilon \ \ \ \ \ \ \ \ \ \ \ \  \ \  \ \text {(3.7)}
	$$
	holds for any  $k \in A$. On the other hand, for any $ k \in B$  we have
	$$
	d\left(T^{k} x, T^{\sigma_{n,m}(k)} y\right) \leq \operatorname{diam}(X)=M .   \ \ \ \ \ \ \ \ \ \ \ \  \ \  \ \text {(3.8)}
	$$
	Hence, we derive that
	$$
	\begin{aligned}
		\inf _{\sigma \in S_{n,m}} \frac{1}{n-m} \sum_{k=m}^{n-1} d\left(T^{k} x, T^{\sigma(k)} y\right) & \leq \frac{1}{n-m} \sum_{k=m}^{n-1}  d\left(T^{k} x, T^{\sigma_{n,m}(k)} y\right) \\
		& =\frac{1}{n-m}  \sum_{k \in A} d\left(T^{k} x, T^{\sigma_{n}(k)} y\right)+\frac{1}{n-m}  \sum_{k \in B} d\left(T^{k} x, T^{\sigma_{n}(k)} y\right) \\
		& \leq \frac{\varepsilon}{n-m} \#(A)+\frac{M}{n-m} \#(B) \ \ \ \ \ \ \text{(according to({3.7}),({3.8}))} \\
		& \leq \frac{\varepsilon}{n-m} \sum_{s=1}^{s_{0}}\left(a_{s}(n-m) +1\right)+\frac{M}{n-m}\left((n-m)-\sum_{\delta=1}^{s_{0}}\left(a_{s}-\delta\right) (n-m)\right).
	\end{aligned}
	$$
	where the last inequality comes from (3.5) and (3.6).\\
	Let  $n-m \rightarrow+\infty $, then we have
	$$
	\overline{BF}(x, y)=\limsup _{n-m \rightarrow+\infty} \inf _{\sigma \in S_{n,m}} \frac{1}{n-m} \sum_{k=m}^{n-1} d\left(T^{k} x, T^{\sigma(k)} y\right) \leq \varepsilon \sum_{s=1}^{s_{0}} a_{s}+M\left(1-\sum_{s=1}^{s_{0}}\left(a_{s}-\delta\right)\right) .$$
	Let  $\delta \rightarrow 0 $, and then we deduce that
	$$
	\overline{BF}(x, y) \leq \varepsilon \sum_{s=1}^{s_{0}} a_{s}+M\left(1-\sum_{s=1}^{s_{0}} a_{s}\right) .
	$$
	The proof  is completed.\\
	With respect to $ \underline{BF}(x, y)$, we have the similar proposition.
	
	{\bf Proposition 3.4} Let  $(X, T)$  be a topological dynamical system, and  $\left\{U_{s}\right\}_{s=1}^{s_{0}} $ be mutually disjoint subsets of $ X$. Given  $x, y \in X$. If the following two conditions hold:
	
	(1) There is  $\varepsilon>0$  such that diam  $\left(U_{s}\right) \leq \varepsilon$  holds for any  $s \in\left\{1,2, \cdots, s_{0}\right\} $;
	
	(2) There is a non-negative integer  interval column  $\left\{[m_r,n_r]\right\}_{r=1}^{\infty}$  of  $\mathbb{N}^{+}$ such that for any  $s \in\left\{1,2, \cdots, s_{0}\right\} $, the following inequalities
	$$
	\lim _{r \rightarrow+\infty} \frac{1}{n_{r}-m_{r}} \sum_{k=m_{r}}^{n_{r}-1} \chi_{U_{g}}\left(T^{k} x\right) \geq a_{s}
	$$
	and
	$$
	\lim _{r \rightarrow+\infty} \frac{1}{n_{r}-m_{r}} \sum_{k=m_{r}}^{n_{r}-1} \chi_{U_{s}}\left(T^{k} y\right) \geq a_{s} .
	$$
	hold for some  $a_{s} \geq 0 .$
	Then, we have
	$$
	\underline{BF}(x, y) \leq \varepsilon \sum_{s=1}^{s_{0}} a_{s}+M\left(1-\sum_{s=1}^{s_{0}} a_{s}\right),
	$$
	where  $M=\operatorname{diam}(X) .$

\section{Uniformly generic points with $ BF(x, y)$}
In this secition, we will prove the Theorem 4.2,  in order to complete the proof process, we assume the contrary that there are uniformly generic points  $x, y \in X $ with $ BF(x, y)$  does not exist, which implies that  $\alpha=\overline{BF}(x, y)-\underline{BF}(x, y)>0 .$ And furthermore we have to estimate  contradiction with the assumption. Then we estimate the upper bound of  $\overline {BF}(x, y)$  and the lower bound of  $\underline{BF}(x, y)$, from which we deduce that  $\overline {BF}(x, y) -\underline{BF}(x, y) \leq \frac{\alpha}{2}$. So  $BF(x, y)$  exists when  $x, y \in X$  are uniformly generic points. In the Proof of Theorem 4.2, we need the following lemma which is a direct corollary of Birkhoff-Von Neumann Theorem [22].

{\bf Lemma 4.1} Let  $X$  be a metric space and  $m, n ,l\in \mathbb{N}^{+}$. If  $\left\{x_{i}\right\}_{i=m}^{n-1},\left\{y_{i}\right\}_{i=m}^{n-1},\left\{\overline{x}_{j}\right\}_{j=m}^{l(n-m)+m-1}$  and  $\left\{\overline{y}_{j}\right\}_{j=m}^{l(n-m)+m-1}$  are subsequences of  $X $ and there is  $\sigma_{l(n-m)+m-1 ,m} \in S_{l(n-m)+m-1 ,m}$  such that
$$
\#\left(\left\{j: \overline{x}_{j}=x_{i}\right\}\right)=\#\left(\left\{j: \overline{y}_{j}=y_{i}\right\}\right)=l
$$
holds for any  $i \in\{m,m+1, \cdots, n-1\}$ , then we have
$$
\inf _{\sigma \in S_{l(n-m)+m-1 ,m}} \sum_{j=m}^{l(n-m)+m-1} d\left(\overline{x}_{j}, \overline{y}_{\sigma(j)}\right)=l \cdot \inf _{\sigma \in S_{n,m}} \sum_{i=m}^{n-1} d\left(x_{i}, y_{\sigma(i)}\right) .
$$

{\bf Theorem 4.2 }Let  $(X, T)$  be a topological dynamical system. If  $x, y \in X$  are uniformly generic points, then  $BF(x, y)$ exists.

{\bf Proof } Let  $x, y \in X $ be uniformly generic points of  $(X, T) $. Put
$$
\alpha=\overline{BF}(x, y)-\underline{BF}(x, y)
$$
We assume that  $BF(x, y)$  does not exist, then  $\alpha>0 .$
Since $ x, y $ are uniformly  generic points, there are  $\mu_{x}, \mu_{y} \in M(X, T) $ such that
$$
\mu_{x}=\lim _{n-m \rightarrow+\infty} \frac{1}{n-m} \sum_{k=m}^{n-1} \delta_{T^{k} x} \text { and } \mu_{y}=\lim _{n-m \rightarrow+\infty} \frac{1}{n-m} \sum_{k=m}^{n-1}  \delta_{T^{k} y} .
$$
Let  $\varepsilon=\frac{\alpha}{8+16 M}$, where  $M=\operatorname{diam}(X)$. According to Lemma 2.3, there exist finite mutually disjoint open sets  $\left\{\Lambda_{r}\right\}_{r=1}^{r_{0}} $ of  $X$  such that
$$
\sum_{r=1}^{r_{0}} \mu_{x}\left(\Lambda_{r}\right) \geq 1-\varepsilon,  \operatorname{diam}\left(\Lambda_{r}\right) \leq \varepsilon, \forall r=1,2, \cdots, r_{0} . \ \ \ \ \ \ \ \text{(4.1)}
$$
Similarly, there exist finite mutually disjoint open sets  $\left\{V_{s}\right\}_{\delta=1}^{\delta_{0}} $ of  $X $ such that
$$
\sum_{s=1}^{s_{0}} \mu_{y}\left(V_{s}\right) \geq 1-\varepsilon, \operatorname{diam}\left(V_{s}\right) \leq \varepsilon, \forall s=1,2, \cdots, s_{0} .\ \ \ \ \ \ \ \text{(4.2)}
$$
Without loss of generality, we can assume that $ \mu_{x}\left(\Lambda_{r}\right)>0$  and $ \mu_{y}\left(V_{s}\right)>0$  for any $ r \in\left\{1,2, \cdots, r_{0}\right\}$,  $s \in\left\{1,2, \cdots, s_{0}\right\} $. Let  $ \Lambda_{r_{0}+1}=X \backslash \bigcup_{r=1}^{r_{0}} \Lambda_{r}$  and $ V_{s_{0}+1}=X \backslash \bigcup_{s=1}^{s_{0}} V_{s} $. We select sequences  $\left\{x_{r}\right\}_{r=1}^{r_{0}+1}$  and  $\left\{y_{s}\right\}_{s=1}^{s_{0}+1}$  of  $X$  such that  $x_{r} \in \Lambda_{r} $ and  $y_{s} \in V_{s} $ for any $ r \in\left\{1,2, \cdots, r_{0}+1\right\}, s \in\left\{1,2, \cdots, s_{0}+1\right\} .$ 

Let  $\left\{\hat{x}_{k}\right\}_{k=m}^{n-1}$  and  $\left\{\hat{y}_{k}\right\}_{k=m}^{n-1}$  be sequences such that  $\hat{x}_{k}=x_{r} $ if  $T^{k} x \in \Lambda_{r} $ and  $\hat{y}_{k}=y_{s}$  if  $T^{k} y \in V_{s}$. According to Proposition 5.1, we have
$$
\begin{aligned}
	& \inf _{\sigma \in S_{n,m}} \frac{1}{n-m} \sum_{k=m}^{n-1} d\left(T^{k} x, T^{\sigma(k)} y\right) \\
	\leq & \inf _{\sigma \in S_{n,m}} \frac{1}{n-m} \sum_{k=m}^{n-1} d\left(T^{k} x, \hat{x}_{\sigma(k)}\right)+\inf _{\sigma \in S_{n,m}} \frac{1}{n-m} \sum_{k=m}^{n-1} d\left(\hat{x}_{\sigma(k)}, \hat{y}_{k}\right)+\inf _{\sigma \in S_{n,m}} \frac{1}{n-m} \sum_{k=m}^{n-1} d\left(\hat{y}_{k}, T^{\sigma(k)} y\right) \\
	= & \inf _{\sigma \in S_{n,m}} \frac{1}{n-m} \sum_{k=m}^{n-1} d\left(\hat{x}_{k}, T^{\sigma(k)} x\right)+\inf _{\sigma \in S_{n,m}} \frac{1}{n-m} \sum_{k=m}^{n-1} d\left(\hat{x}_{k}, \hat{y}_{\sigma(k)}\right)+\inf _{\sigma \in S_{n,m}} \frac{1}{n-m} \sum_{k=m}^{n-1} d\left(\hat{y}_{k}, T^{\sigma(k)} y\right) . \ \ \ \ \ \ \ \text{(4.3)}
\end{aligned}
$$
Similarly, we have
$$
\begin{aligned}
	& \inf _{\sigma \in S_{n,m}} \frac{1}{n-m} \sum_{k=m}^{n-1} d\left(\hat{x}_{k}, \hat{y}_{\sigma(k)}\right) \\
	\leq & \inf _{\sigma \in S_{n,m}} \frac{1}{n-m} \sum_{k=m}^{n-1} d\left(\hat{x}_{k}, T^{\sigma(k)} x\right)+\inf _{\sigma \in S_{n,m}} \frac{1}{n-m} \sum_{k=m}^{n-1} d\left(T^{\sigma(k)} x, T^{k} y\right)+\inf _{\sigma \in S_{n,m}} \frac{1}{n-m} \sum_{k=m}^{n-1} d\left(T^{k} y, \hat{y}_{\sigma(k)}\right) \\
	= & \inf _{\sigma \in S_{n,m}} \frac{1}{n-m} \sum_{k=m}^{n-1} d\left(\hat{x}_{k}, T^{\sigma(k)} x\right)+\inf _{\sigma \in S_{n,m}} \frac{1}{n-m} \sum_{k=m}^{n-1} d\left(T^{k} x, T^{\sigma(k)} y\right)+\inf _{\sigma \in S_{n,m}} \frac{1}{n-m} \sum_{k=m}^{n-1} d\left(\hat{y}_{k}, T^{\sigma(k)} y\right) .
\end{aligned}
$$
Thus, we deduce that
$$
\begin{aligned}
	& \inf _{\sigma \in S_{n,m}} \frac{1}{n-m} \sum_{k=m}^{n-1} d\left(T^{k} x, T^{\sigma(k)} y\right) \\
	\geq & \inf _{\sigma \in S_{n,m}} \frac{1}{n-m} \sum_{k=m}^{n-1} d\left(\hat{x}_{k}, \hat{y}_{\sigma(k)}\right)-\inf _{\sigma \in S_{n,m}} \frac{1}{n-m} \sum_{k=m}^{n-1} d\left(\hat{x}_{k}, T^{\sigma(k)} x\right)-\inf _{\sigma \in S_{n,m}} \frac{1}{n-m} \sum_{k=m}^{n-1} d\left(\hat{y}_{k}, T^{\sigma(k)} y\right) .\ \ \ \ \ \ \ \text{(4.4)}
\end{aligned}
$$
Our idea is to estimate the bounds of  $\inf _{\sigma \in S_{n,m}} \frac{1}{n-m} \sum_{k=m}^{n-1}  d\left(T^{k} x, T^{\sigma(k)} y\right)$.\\ 
Inequalities (4.3) and (4.4) show that we only need to estimate the bounds of  $\inf _{\sigma \in S_{n,m}} \frac{1}{n-m} \sum_{k=m}^{n-1} d\left(\hat{x}_{k}, \hat{y}_{\sigma(k)}\right)$,  $\inf _{\sigma \in S_{n,m}} \frac{1}{n-m} \sum_{k=m}^{n-1} d\left(\hat{x}_{k}, T^{\sigma(k)} x\right) $ and  $\inf _{\sigma \in S_{n,m}} \frac{1}{n-m} \sum_{k=m}^{n-1}  d\left(\hat{y}_{k}, T^{\sigma(k)} y\right)$. 

In the following, Lemma 4.3 shows the upper bounds of  $\inf _{\sigma \in S_{n,m}} \frac{1}{n-m} \sum_{k=m}^{n-1}  d\left(\hat{x}_{k}, T^{\sigma(k)} x\right) $ and  $\inf_{\sigma \in S_{n,m}} \frac{1}{n-m} \sum_{k=m}^{n-1}  d\left(\hat{y}_{k}, T^{\sigma(k)} y\right)$, Lemma 6.3 shows the lower bound of  $\inf _{\sigma \in S_{n,m}} \frac{1}{n-m} \sum_{k=m}^{n-1}  d\left(\hat{x}_{k}, \hat{y}_{\sigma(k)}\right)$  and Lemma 4.4 shows the upper bound of  $\inf _{\sigma \in S_{n,m}} \frac{1}{n-m} \sum_{k=m}^{n-1}  d\left(\hat{x}_{k}, \hat{y}_{\sigma(k)}\right)$. 

Given $\beta>0$. According to Lemma 2.1, for all sufficiently large  $n-m \in \mathbb{N}^{+}$, then we have
$$
\frac{1}{n-m} \sum_{k=m}^{n-1} \chi_{\Lambda_{r}}\left(T^{k} x\right) \geq \mu_{x}\left(\Lambda_{r}\right)-\beta$$ \text { and } $$\frac{1}{n-m} \sum_{k=m}^{n-1} \chi_{V_{x}}\left(T^{k} y\right) \geq \mu_{y}\left(V_{s}\right)-\beta       \ \ \ \ \ \ \ \text{(4.5)}
$$
hold for any  $r \in\left\{1,2, \cdots, r_{0}\right\} $ and $ s \in\left\{1,2, \cdots, s_{0}\right\} .$

{\bf Lemma 4.3}
 For all sufficiently large  $n-m \in \mathbb{N}^{+} $, we have
$$
\inf _{\sigma \in S_{n,m}} \frac{1}{n-m} \sum_{k=m}^{n-1} d\left(\hat{x}_{k}, T^{\sigma(k)} x\right) \leq \varepsilon+M \varepsilon+M r_{0}\beta$$  \text { and } $$\inf _{\sigma \in S_{n,m}} \frac{1}{n-m} \sum_{k=m}^{n-1} d\left(\hat{y}_{k}, T^{\sigma(k)} y\right) \leq \varepsilon+M \varepsilon+M s_{0} \beta .
$$
{\bf Proof}  According to (6.1), we have $ d\left(\hat{x}_{k}, T^{k} x\right) \leq \varepsilon $ if $ T^{k} x \in \bigcup_{r=1}^{r_{0}} \Lambda_{r}$. On the other hand, we
can estimate  $d\left(\hat{x}_{k}, T^{k} x\right) \leq M $ for  $\operatorname{diam}(X)=M$. Then we have
$$
\begin{aligned}
	& \inf _{\sigma \in S_{n,m}} \frac{1}{n-m} \sum_{k=m}^{n-1} d\left(\hat{x}_{k}, T^{\sigma(k)} x\right) \leq \frac{1}{n-m} \sum_{k=m}^{n-1} d\left(\hat{x}_{k}, T^{k} x\right) \\
	\leq & \varepsilon \cdot \frac{\#\left(\left\{k \in \mathbb{N}^{+}: T^{k} x \in \bigcup_{r=1}^{r_{0}} \Lambda_{r}, k \in [m,n-1]\right\}\right)}{n-m}+M \cdot \frac{\#\left(\left\{k \in \mathbb{N}^{+}: T^{k} x \notin \bigcup_{r=1}^{r_{0}} \Lambda_{r}, k \in [m,n-1]\right\}\right)}{n-m} \\
	\leq & \varepsilon+M\left(1-\frac{\#\left(\left\{k \in \mathbb{N}^{+}: T^{k} x \in \bigcup_{r=1}^{r_{0}} \Lambda_{r}, k \in [m,n-1]\right\}\right)}{n-m}\right) . \ \ \ \ \ \ \ \text{(4.6)}
\end{aligned}
$$
Since  $\left\{\Lambda_{r}\right\}_{r=1}^{r_{0}}$  are mutually disjoint, for sufficiently large  $n-m \in \mathbb{N}^{+}$ we have
$$
\begin{aligned}
	\frac{\#\left(\left\{k \in \mathbb{N}^{+}: T^{k} x \in \bigcup_{r=1}^{r_{0}} \Lambda_{r}, k\in [m,n-1]\right\}\right)}{n-m} & =\sum_{r=1}^{r_{0}} \frac{\#\left(\left\{k \in \mathbb{N}^{+}: T^{k} x \in \Lambda_{r}, k\in [m,n-1]\right\}\right)}{n-m} \\
	& =\sum_{r=1}^{r_{0}} (\frac{1}{n-m} \sum_{k=m}^{n-1} \chi_{\Lambda_{r}}\left(T^{k} x\right)) \\
	& \geq \sum_{r=1}^{r_{0}}\left(\mu_{x}\left(\Lambda_{r}\right)-\beta\right)(\mathrm{b} y(\text { (4.5 })) \\
	& =\sum_{r=1}^{r_{0}} \mu_{x}\left(\Lambda_{r}\right)-r_{0} \beta .         \ \ \ \ \ \ \ \ \ \ \ \ \ \ \ \ \ \ \ \ \ \text{(4.7)}
\end{aligned}
$$
Combining (4.6) with (4.7), we derive that
$$
\begin{aligned}
	\inf _{\sigma \in S_{n,m}} \frac{1}{n-m} \sum_{k=m}^{n-1} d\left(\hat{x}_{k}, T^{\sigma(k)} x\right) & \leq \varepsilon+M\left(1-\frac{\#\left(\left\{k \in \mathbb{N}^{+}: T^{k} x \in \bigcup_{r=1}^{r_{0}} \Lambda_{r}, k\in [m,n-1]\right\}\right)}{n-m}\right) \\
	& \leq \varepsilon+M\left(1-\sum_{r=1}^{r_{0}} \mu_{x}\left(\Lambda_{r}\right)+r_{0} \beta\right) \\
	& \leq \varepsilon+M \varepsilon+M r_{0} \beta, \ \ \ \ \ \ \ \ \ \ \ \ \ \ \text{(4.8)}
\end{aligned}
$$
where the last inequality comes from (4.1). Similarly, we have
$$
	\inf _{\sigma \in S_{n,m}} \frac{1}{n-m} \sum_{k=m}^{n-1}  d\left(\hat{y}_{k}, T^{\sigma(k)} y\right) \leq \varepsilon+M \varepsilon+M s_{0} \beta . \ \ \ \ \ \ \ \ \ \ \ \ \ \ \text{(4.9)}
$$
This finishes the proof of Lemma 4.3.

To estimate the bounds of  $	\inf _{\sigma \in S_{n,m}} \frac{1}{n-m} \sum_{k=m}^{n-1}   d\left(\hat{x}_{k}, \hat{y}_{\sigma(k)}\right) $, we introduce some notations. Let  $a_{1}=   \min \left\{\mu_{x}\left(\Lambda_{r}\right), r=1,2, \cdots, r_{0}\right\}$, $a_{2}=\min \left\{\mu_{y}\left(V_{s}\right), s=1,2, \cdots, s_{0}\right\} $ and  $a=\min \left\{\frac{\varepsilon}{r_{0}+s_{0}}, a_{1}, a_{2}\right\} .$ Then for any  $r \in\left\{1,2, \cdots, r_{0}\right\}$  and any  $s \in\left\{1,2, \cdots, s_{0}\right\}$ , there are $ n_{r}, m_{s} \in \mathbb{N}^{+} $such that
$$
a n_{r} \leq \mu_{x}\left(\Lambda_{r}\right)<a\left(n_{r}+1\right) \ \ \ \ \ \ \ \ \ \ \ \ \ \ \text{(4.10)}
$$
and
$$
a m_{s} \leq \mu_{y}\left(V_{s}\right)<a\left(m_{s}+1\right) . \ \ \ \ \ \ \ \ \ \ \ \ \ \ \text{(4.11)}
$$
Let  $K=\min \left\{\sum_{r=1}^{r_{0}} n_{r}, \sum_{s=1}^{s_{0}} m_{s}\right\} $. Without loss of generality, we can assume  $K=\sum_{r=1}^{r_{0}} n_{r} \leq \sum_{s=1}^{s_{0}} m_{s}$ . Then there is  $s_{0}^{*} \leq s_{0}$  and integer sequence  $\left\{m_{s}^{*}\right\}_{s=1}^{s_{0}^{*}} $ such that
$$
1 \leq m_{s}^{*} \leq m_{s}, \forall s=1,2, \cdots, s_{0}^{*} \ \ \ \ \ \ \ \ \ \ \ \ \ \ \text{(4.12)}
$$
and
$$
K=\sum_{s=1}^{s_{0}^{*}} m_{s}^{*}     .                   \ \ \ \ \ \ \ \ \ \ \ \ \ \ \text{(4.13)}
$$
According to (4.1) we derive that
$$
\begin{aligned}
	1-K a & \leq \varepsilon+\sum_{r=1}^{r_{0}} \mu_{x}\left(\Lambda_{r}\right)-a \sum_{r=1}^{r_{0}} n_{r} \\
	& =\varepsilon+\sum_{r=1}^{r_{0}}\left(\mu_{x}\left(\Lambda_{r}\right)-a n_{r}\right)                            \ \ \ \  \ \ \ \ \text{(4.14)}\\
	& <\varepsilon+r_{0} a  \ \ \ (\mathrm{according\  to}(4.10)) \\
	& \leq 2 \varepsilon,
\end{aligned}
$$
where the last inequality comes from  $a \leq \frac{\varepsilon}{r_{0}+\varepsilon_{0}} .$
Construct the sequences  $\left\{\overline{x}_{i}\right\}_{i=m}^{m+K-1}$  and $ \left\{\overline{y}_{i}\right\}_{i=m}^{m+K-1} $ as follows:
$$
\overline{x}_{i}=\left\{\begin{array}{ll}
	x_{1}, &  m \leq i \leq m+n_{1}-1 \\
	x_{r+1}, & \sum_{j=1}^{r} n_{j}+m-1<i \leq \sum_{j=1}^{r+1} n_{j}+m-1, r=1,2, \cdots, r_{0}-1,
\end{array}\right.
$$
$$
\overline{y}_{i}=\left\{\begin{array}{ll}
	y_{1}, & m \leq i\leq m_{1}^{*}+m-1 \\
	y_{s+1}, & \sum_{j=1}^{s} m_{j}^{*}+m-1<i \leq \sum_{j=1}^{s+1} m_{j}^{*}+m-1, s=1,2, \cdots, s_{0}^{*}-1 .
\end{array}\right.
$$
In the following we will use the distance between sequences  $\left\{\overline{x}_{i}\right\}_{i=m}^{m+K-1}$  and  $\left\{\overline{y}_{i}\right\}_{i=m}^{m+K-1} $ to estimate the bounds of  $\inf _{\sigma \in S_{n}} \frac{1}{n-m} \sum_{k= m}^{n-1} d\left(\hat{x}_{k}, \hat{y}_{\sigma(k)}\right).$

{\bf Lemma 4.4} For all sufficiently large  $n-m \in \mathbb{N}^{+}$, we have
$$
\inf _{\sigma \in S_{n,m}} \frac{1}{n-m} \sum_{k=m}^{n-1} d\left(\hat{x}_{k}, \hat{y}_{\sigma(k)}\right) \geq(a-\beta) \cdot \inf _{\sigma \in S_{K+m,m}} \sum_{k=m}^{K+m-1} d\left(\overline{x}_{k}, \overline{y}_{\sigma(k)}\right)-\frac{M K}{n-m}-2 M \varepsilon-M K \beta .
$$
{\bf Proof} Let  $l_{n,m}$  be the minimal integer such that  $l_{n,m} \geq(a-\beta)(n-m) $. Then according to (4.14) we have
$$
\frac{(n-m)-K l_{n,m}}{(n-m)} \leq 1-K(a-\beta) \leq 2 \varepsilon+K \beta  .            \ \ \ \  \ \ \ \ \text{(4.15)}
$$

{\bf Claim 1}
  For all sufficiently large  $n-m \in \mathbb{N}^{+}$ , we have
$$
\#\left(\left\{k \mid \hat{x}_{k}=x_{r}, k \in[m,n-1] \right\}\right) \geq l_{n,m} n_{r}$$ \text { and } $$\#\left(\left\{k \mid \hat{y}_{k}=y_{s},k \in[m,n-1]  \right\}\right) \geq l_{n,m} m_{s}^{*}
$$
hold for any  $r \in\left\{1,2, \cdots, r_{0}\right\} $ and any $ s \in\left\{1,2, \cdots, s_{0}^{*}\right\} .$

{\bf Proof of claim 1}  Combining (4.5) with (4.10), for sufficiently large $ n-m \in \mathbb{N}+ $, we have
$$
\frac{1}{n-m} \sum_{k=m}^{n-1} \chi_{\Lambda_{r}}\left(T^{k} x\right) \geq \mu_{x}\left(\Lambda_{r}\right)-\beta \geq a n_{r}-\beta
$$
holds for any  $r \in\left\{1,2, \cdots, r_{0}\right\} .$
If  $n_{r}=1$ , we have
$$
\sum_{k=m}^{n-1} \chi_{\Lambda_{r}}\left(T^{k} x\right) \geq a (n-m) n_{r}-\beta (n-m)=(a-\beta)(n-m)
$$
Since  $\sum_{k=m}^{n-1} \chi_{\Lambda_{r}}\left(T^{k} x\right)$  is integer and  $l_{n,m}$  is the minimal integer such that  $l_{n,m} \geq(a-\beta)(n-m)$, we derive that
$$
\#\left(\left\{k \mid \hat{x}_{k}=x_{r}, k\in[m,n-1]\right\}\right)=\sum_{k=m}^{n-1} \chi_{\Lambda_{r}}\left(T^{k} x\right) \geq l_{n,m}=l_{n,m} n_{r}
$$
If  $n_{r}>1$, for sufficiently large  $n-m \in \mathbb{N}+$  we have
$$
\beta (n-m)\left(n_{r}-1\right)-n_{r}>0
$$
Then we deduce that
$$
\begin{aligned}
	\sum_{k=m}^{n-1} \chi_{\Lambda_{r}}\left(T^{k} x\right) & \geq a (n-m) n_{r}-\beta (n-m) \\
	& =l_{n,m} n_{r}-l_{n,m} n_{r}+a (n-m) n_{r}-\beta (n-m) n_{r}+\beta (n-m)\left(n_{r}-1\right) \\
	& =l_{n,m} n_{r}+\beta (n-m)\left(n_{r}-1\right)-\left(l_{n,m}-(a-\beta) (n-m)\right) n_{r} \\
	& \geq l_{n,m} n_{r}+\beta (n-m)\left(n_{r}-1\right)-n_{r} \\
	& >l_{n,m} n_{r}.
\end{aligned}
$$
and 
$$
\#\left(\left\{k \mid \hat{x}_{k}=x_{r}, k \in [m,n-1]\right\}\right)=\sum_{k=m}^{n-1} \chi_{\Lambda_{r}}\left(T^{k} x\right)>l_{n,m} n_{r}.$$
Similarly, we can prove that for sufficiently large  $n-m \in \mathbb{N}^{+}$ , we have
$$
\#\left(\left\{k \mid \hat{y}_{k}=y_{s}, k\in [m,n-1]\right\}\right) \geq l_{n,m} m_{s}^{*}
$$
holds for any  $s \in\left\{1,2, \cdots, s_{0}^{*}\right\} .$

According to Claim 1, for sufficiently large $ n-m \in \mathbb{N}^{+}$ there is a set  $A \subset [m,n-1] $ such that
$$
\#(A)=K l_{n,m}    \ \ \ \  \ \ \ \ \text{(4.16)}
$$
and
$$
\#\left(\left\{k \mid \hat{x}_{k}=x_{r}, k \in A\right\}\right)=n_{r}l_{n,m} \ \ \ \  \ \ \ \ \text{(4.17)}
$$
holds for any  $r \in\left\{1,2, \cdots, r_{0}\right\} .$

Let  $\sigma_{1} \in S_{n,m}$  such that
$$
\frac{1}{n-m} \sum_{k=m}^{n-1} d\left(\hat{x}_{k}, \hat{y}_{\sigma_{1}(k)}\right)=\inf _{\sigma \in S_{n,m}} \frac{1}{n-m} \sum_{k=m}^{n-1} d\left(\hat{x}_{k}, \hat{y}_{\sigma(k)}\right) .
$$
For any $ s \in\left\{1,2, \cdots, s_{0}^{*}\right\}$, we denote
$$
b_{s}=\#\left(\left\{k \in [m,n-1] \mid \hat{y}_{\sigma_{1}(k)}=y_{s}, k \notin A\right\}\right)
$$
Then according to (4.16), we have
$$
\sum_{s=1}^{s_{0}^{*}} b_{s} \leq n-m-\#(A)=n-m-K l_{n,m}          \ \ \ \  \ \ \ \ \text{(4.18)}
$$
Hence,for any  $s \in\left\{1,2, \cdots, s_{0}^{*}\right\} $ and sufficiently large $ n-m \in \mathbb{N}^{+}$ , we have
$$
\begin{aligned}
	\frac{\#\left(\left\{k \in \mathbb{N}^{+} \mid \hat{y}_{\sigma_{1}(k)}=y_{s}, k \in A\right\}\right)+b_{s}}{n-m} &=\frac{\#\left(\left\{k \in \mathbb{N}^{+} \mid \hat{y}_{\sigma_{1}(k)}=y_{s}, k\in [m,n-1]\right\}\right)}{n-m} \\
	 &=\frac{\#({k \in \mathbb{N}^{+} \mid \hat{y}_{k}=y_{s}, k \in [m,n-1]})}{n-m}\\
	 &=\frac{1}{n-m} \sum_{k=m}^{n-1} \chi_{V_{s}}(T^{k} y) \\
	 &\geq \mu_{y}\left(V_{s}\right)-\beta \ \ \ \ \ \ \ (\text{by 4.5} )\\
	 &\geq m_{s} a-\beta \ \ \ \ \ \ \  (\text{by (4.11)})) \\
	 &\geq m_{s}^{*}(a-\beta) \\
\end{aligned}
$$
where the last inequality follows from (4.12). Therefore,
$$
\frac{\#\left(\left\{k \in \mathbb{N}^{+} \mid \hat{y}_{\sigma_{1}(k)}=y_{s}, k \in A\right\}\right)}{n-m} \geq m_{s}^{*}(a-\beta)-\frac{b_{s}}{n-m} . \ \ \  \ \ \ \ \text{(4.19)}
$$
According to Claim 1 and (6.19), for sufficiently large $ n-m \in \mathbb{N}^{+}$, there is a  $\sigma_{2} \in S_{n,m}$  such that
$$
\#\left(\left\{k \mid \hat{y}_{\sigma_{2}(k)}=y_{s}, k \in A\right\}\right)=l_{n,m} m_{s}^{*}    \ \ \  \ \ \ \ \text{(4.20)}
$$
and
$$
\frac{\#\left(\left\{k \in \mathbb{N}^{+} \mid \sigma_{2}(k)=\sigma_{1}(k), \hat{y}_{\sigma_{1}(k)}=y_{s}, k \in A\right\}\right)}{n-m} \geq m_{s}^{*}(a-\beta)-\frac{b_{s}}{n-m}$$
hold for any  $s \in\left\{1,2, \cdots, s_{0}^{*}\right\}$.

 Since  $\#(A)=K l_{n,m}$, we deduce that
$$
\begin{aligned}
	\#\left(\left\{k \in \mathbb{N}^{+} \mid \sigma_{2}(k) \neq \sigma_{1}(k), k \in A\right\}\right) & =\#(A)-\#\left(\left\{k \in \mathbb{N}^{+} \mid \sigma_{2}(k)=\sigma_{1}(k), k \in A\right\}\right) \\
	& \leq K l_{n,m}-\sum_{s=1}^{s_{0}^{*}} \#\left(\left\{k \in \mathbb{N}^{+} \mid \sigma_{2}(k)=\sigma_{1}(k), \hat{y}_{\sigma_{1}(k)}=y_{s}, k \in A\right\}\right) \\
	& \leq K l_{n,m}-(n-m) \sum_{s=1}^{s_{0}^{*}}\left(m_{s}^{*}(a-\beta)-\frac{b_{s}}{n-m}\right) \\
	& =K l_{n,m}-\sum_{s=1}^{s_{0}^{*}} m_{s}^{*}(a-\beta) (n-m)+\sum_{s=1}^{s_{0}^{*}} b_{s} \\
	& =K\left(l_{n,m}-(a-\beta) (n-m)\right)+\sum_{s=1}^{s_{0}^{*}} b_{s} \text { (according to (4.13) ) } \\
	& \leq K+n-m-K l_{n,m},
\end{aligned}
$$
where the last inequality comes from (4.18) and the fact that $ l_{n,m}$  is the minimal integer such that  $l_{n,m} \geq(a-\beta)(n-m)$. Then we derive that
$$
\begin{aligned}
	\sum_{k \in A} d\left(\hat{x}_{k}, \hat{y}_{\sigma_{2}(k)}\right) & =\sum_{k \in A, \sigma_{1}(k)=\sigma_{2}(k)} d\left(\hat{x}_{k}, \hat{y}_{\sigma_{2}(k)}\right)+\sum_{k \in A, \sigma_{1}(k) \neq \sigma_{2}(k)} d\left(\hat{x}_{k}, \hat{y}_{\sigma_{2}(k)}\right) \\
	& \leq \sum_{k \in A} d\left(\hat{x}_{k}, \hat{y}_{\sigma_{1}(k)}\right)+M \cdot \#\left(\left\{k \in \mathbb{N}^{+} \mid \sigma_{2}(k) \neq \sigma_{1}(k), k \in A\right\}\right) \\
	& \leq \sum_{k \in A} d\left(\hat{x}_{k}, \hat{y}_{\sigma_{1}(k)}\right)+M\left(K+n-m-K l_{n,m}\right) .
\end{aligned}
$$
This implies that
$$
\begin{aligned}
	\frac{1}{n-m} \sum_{k \in A} d\left(\hat{x}_{k}, \hat{y}_{\sigma_{2}(k)}\right) & \leq \frac{1}{n-m} \sum_{k \in A} d\left(\hat{x}_{k}, \hat{y}_{\sigma_{1}(k)}\right)+\frac{M K}{n-m}+\frac{M\left(n-m-K l_{n,m}\right)}{n-m} \\
	& \leq \frac{1}{n-m} \sum_{k \in A} d\left(\hat{x}_{k}, \hat{y}_{\sigma_{1}(k)}\right)+\frac{M K}{n-m}+2 M \varepsilon+M K \beta,
\end{aligned}
$$
where the last inequality comes from (4.15).

Let  $B=\left\{\sigma_{2}(k) \mid k \in A\right\}$  and $ P(A, B)=\left\{\sigma \in S_{n,m} \mid \sigma(k) \in B, \forall k \in A\right\}.$  Then  $\sigma_{2} \in P(A, B) $ and we have
$$
\begin{aligned}
	\inf _{\sigma \in S_{n,m}} \frac{1}{n-m} \sum_{k=m}^{n-1} d\left(\hat{x}_{k}, \hat{y}_{\sigma(k)}\right) & =\frac{1}{n-m} \sum_{k=m}^{n-1} d\left(\hat{x}_{k}, \hat{y}_{\sigma_{1}(k)}\right) \geq \frac{1}{n-m} \sum_{k \in A} d\left(\hat{x}_{k}, \hat{y}_{\sigma_{1}(k)}\right) \\
	& \geq \frac{1}{n-m} \sum_{k \in A} d\left(\hat{x}_{k}, \hat{y}_{\sigma_{2}(k)}\right)-\frac{M K}{n-m}-2 M \varepsilon-M K \beta \\
	& \geq \inf _{\sigma \in P(A, B)} \frac{1}{n-m} \sum_{k \in A} d\left(\hat{x}_{k}, \hat{y}_{\sigma(k)}\right)-\frac{M K}{n-m}-2 M \varepsilon-M K \beta .\ \ \  \ \ \ \ \text{(4.21)}
\end{aligned}
$$
According to (4.17), (4.20) and Lemma 4.1 we can derive that
$$
\inf _{\sigma \in P(A, B)} \sum_{k \in A} d\left(\hat{x}_{k}, \hat{y}_{\sigma(k)}\right)=l_{n,m} \cdot \inf _{\sigma \in S_{K+m,m}} \sum_{k=m}^{K+m-1} d\left(\overline{x}_{k}, \overline{y}_{\sigma(k)}\right).   \ \ \ \ \text{(4.22)}$$

Combining (4.21) with (4.22), we estimate the lower bound of  $\inf _{\sigma \in S_{n,m}} \frac{1}{n-m} \sum_{k=m}^{n-1} d\left(\hat{x}_{k}, \hat{y}_{\sigma(k)}\right) $ as follows:
$$
\begin{aligned}
	\inf _{\sigma \in S_{n,m}} \frac{1}{n-m} \sum_{k=m}^{n-1} d\left(\hat{x}_{k}, \hat{y}_{\sigma(k)}\right) & \geq \frac{l_{n,m}}{n-m} \cdot \inf _{\sigma \in S_{K+m,m}} \sum_{k=m}^{K+m-1} d\left(\overline{x}_{k}, \overline{y}_{\sigma(k)}\right)-\frac{M K}{n-m}-2 M \varepsilon-M K \beta \\
	& \geq(a-\beta) \cdot \inf _{\sigma \in S_{K+m,m}} \sum_{k=m}^{K+m-1} d\left(\overline{x}_{k}, \overline{y}_{\sigma(k)}\right)-\frac{M K}{n-m}-2 M \varepsilon-M K \beta .
\end{aligned}
$$
The proof of Lemma 4.4 is completed.

{\bf Lemma 4.5} For all sufficiently large  $n-m \in \mathbb{N}^{+}$, we have
$$
\inf _{\sigma \in S_{n,m}} \frac{1}{n-m} \sum_{k=m}^{n-1} d\left(\hat{x}_{k}, \hat{y}_{\sigma(k)}\right) \leq(a-\beta) \cdot \inf _{\sigma \in S_{K+m,m}} \sum_{k=m}^{K+m-1} d\left(\overline{x}_{k}, \overline{y}_{\sigma(k)}\right)+\frac{M K}{n-m}+2 M \varepsilon+M K \beta .
$$
{\bf Proof} Let $ \sigma_{3} \in S_{K+m,m}$  such that
$$
\sum_{k=m}^{K+m-1}d\left(\overline{x}_{k}, \overline{y}_{\sigma_{3}(k)}\right)=\inf _{\sigma \in S_{K+m,m}} \sum_{k=m}^{K+m-1} d\left(\overline{x}_{k}, \overline{y}_{\sigma(k)}\right)
$$
For sufficiently large  $n-m \in \mathbb{N}^{+} $, Claim 1 shows that there is a partition  $\left\{A_{i}\right\}_{i=m}^{m+K}$  of  $\{m,m+1, \cdots, n-1\}$  such that
$$
\#\left(A_{i}\right)=l_{n,m} \text { and }\left\{\hat{x}_{k} \mid k \in A_{i}\right\}=\left\{\overline{x}_{i}\right\}
\ \ \text{(4.23)}$$
hold for any  $i \in\{m,m+1, \cdots, m+K-1\} $. 

Similarly, there is a partition  $\left\{B_{i}\right\}_{i=m}^{m+K}$  of $ \{m,m+1, \cdots, n-1 \}$  such that
$$
\#\left(B_{i}\right)=l_{n,m} \text { and }\left\{\hat{y}_{k} \mid k \in B_{i}\right\}=\left\{\overline{y}_{\sigma_{3}(i)}\right\}
\ \ \text{(4.24)}$$
hold for any  $i \in\{m,m+1, \cdots, m+K-1\}.$

According to (4.23) and (4.24), there exists  $\sigma_{n,m} \in S_{n,m}$  such that
$$
B_{i}=\left\{\sigma_{n,m}(k) \mid k \in A_{i}\right\}, \forall i=m,m+1, \cdots, m+K-1.
$$
Thus, for any $ i \in\{m,m+1, \cdots, m+K-1\} $, we have that
$$
\sum_{k \in A_{i}} d\left(\hat{x}_{k}, \hat{y}_{\sigma_{{n,m}}(k)}\right)=l_{n,m} \cdot d\left(\overline{x}_{i}, \overline{y}_{\sigma_{3}(i)}\right)
$$
Hence, we deduce that
$$
\begin{aligned}
	\sum_{k=m}^{n-1} d\left(\hat{x}_{k}, \hat{y}_{\sigma_{n,m}(k)}\right) & =\sum_{i=m}^{m+K} \sum_{k \in A_{i}} d\left(\hat{x}_{k}, \hat{y}_{\sigma_{n,m}(k)}\right) \\
	& =\sum_{i=m}^{K+m-1} l_{n,m} \cdot d\left(\overline{x}_{i}, \overline{y}_{\sigma_{3}(i)}\right)+\sum_{k \in A_{K+m}} d\left(\hat{x}_{k}, \hat{y}_{\sigma_{n,m}(k)}\right) \\
	& \leq l_{n,m} \cdot \sum_{i=m}^{K+m-1} d\left(\overline{x}_{i}, \overline{y}_{\sigma_{3}(i)}\right)+M\left(n-m-K l_{n,m}\right) \\
	& \leq\left(l_{n,m}-1\right) \cdot \sum_{k=m}^{m+K-1} d\left(\overline{x}_{k}, \overline{y}_{\sigma_{3}(k)}\right)+M K+M\left(n-m-K l_{n,m}\right) .
\end{aligned}
$$

Since $l_{n,m}-1<(a-\beta)(n-m)$ , we have
$$
\begin{aligned}
	\frac{1}{n-m} \sum_{k=m}^{n-1} d\left(\hat{x}_{k}, \hat{y}_{\sigma_{n,m}(k)}\right) & \leq(a-\beta) \cdot \sum_{k=m}^{K+m-1} d\left(\overline{x}_{k}, \overline{y}_{\sigma_{3}(k)}\right)+\frac{M K}{n-m}+\frac{M\left(n-m-K l_{n,m}\right)}{n-m} \\
	& \leq(a-\beta) \cdot  \sum_{k=m}^{K+m-1} d\left(\overline{x}_{k}, \overline{y}_{\sigma_{3}(k)}\right)+\frac{M K}{n-m}+2 M \varepsilon+M K \beta. \ \text{(According to 6.15)}
\end{aligned}
$$
 Hence, we can estimate the upper bound of  $\inf _{\sigma \in S_{n,m}} \frac{1}{n-m} \sum_{k=m}^{n-1} d\left(\hat{x}_{k}, \hat{y}_{\sigma(k)}\right) $ as follows:
$$
\begin{aligned}
	\inf _{\sigma \in S_{n,m}} \frac{1}{n-m} \sum_{k=m}^{n-1} d\left(\hat{x}_{k}, \hat{y}_{\sigma(k)}\right) & \leq  \frac{1}{n-m} \sum_{k=m}^{n-1} d\left(\hat{x}_{k}, \hat{y}_{\sigma_{n,m}(k)}\right) \\
	& \leq(a-\beta) \cdot \inf _{\sigma \in S_{K+m,m}} \sum_{k=m}^{K+m-1} d\left(\overline{x}_{k}, \overline{y}_{\sigma(k)}\right)+\frac{M K}{n-m}+2 M \varepsilon+M K \beta .
\end{aligned}
$$
This finishes the proof of Lemma 4.5.

Given sufficiently large  $n-m \in \mathbb{N}^{+}$. According to (6.3), Lemma 4.3 and Lemma 4.5, we deduce that
$$
\begin{aligned}
	& \inf _{\sigma \in S_{n,m}} \frac{1}{n-m} \sum_{k=m}^{n-1} d\left(T^{k} x, T^{\sigma(k)} y\right) \\
	\leq & (a-\beta) \cdot \inf _{\sigma \in S_{K+m,m}} \sum_{k=m}^{K+m-1} d\left(\overline{x}_{k}, \overline{y}_{\sigma(k)}\right)+\frac{M K}{n-m}+2 M \varepsilon+M K \beta+2 \varepsilon+2 M \varepsilon+M\left(r_{0}+s_{0}\right) \beta .
\end{aligned}
$$
On the other hand, according to (4.4), Lemma 4.3 and Lemma 4.4, we derive that
$$
\begin{aligned}
	& \inf _{\sigma \in S_{n,m}} \frac{1}{n-m} \sum_{k=m}^{n-1} d\left(T^{k} x, T^{\sigma(k)} y\right) \\
	\geq & (a-\beta) \cdot \inf _{\sigma \in S_{K+m,m}} \sum_{k=m}^{K+m-1} d\left(\overline{x}_{k}, \overline{y}_{\sigma(k)}\right)-\frac{M K}{n-m}-2 M \varepsilon-M K \beta-2 \varepsilon-2 M \varepsilon-M\left(r_{0}+s_{0}\right) \beta .
\end{aligned}
$$
Let  $n-m \rightarrow+\infty$, and then we deduce that
$$
\overline{BF}(x, y) \leq(a-\beta) \cdot \inf _{\sigma \in S_{K+m,m}} \sum_{k=m}^{K+m-1} d\left(\overline{x}_{k}, \overline{y}_{\sigma(k)}\right)+2 M \varepsilon+M K \beta+2 \varepsilon+2 M \varepsilon+M\left(r_{0}+s_{0}\right) \beta
$$
and
$$
\underline{BF}(x, y) \geq(a-\beta) \cdot \inf _{\sigma \in S_{K+m,m}} \sum_{k=m}^{K+m-1} d\left(\overline{x}_{k}, \overline{y}_{\sigma(k)}\right)-2 M \varepsilon-M K \beta-2 \varepsilon-2 M \varepsilon-M\left(r_{0}+s_{0}\right) \beta .
$$
This shows

$$ \overline{BF}(x, y)-\underline{BF}(x, y) \leq 8 M \varepsilon+2 M K \beta+4 \varepsilon+2 M\left(r_{0}+s_{0}\right) \beta .$$
Let  $\beta \rightarrow 0 $, then we derive that
$$
\overline{BF}(x, y)-\underline{BF}(x, y) \leq 8 M \varepsilon+4 \varepsilon=\frac{\alpha}{2},
$$
which is a contradiction with the assumption. Hence,  $BF(x, y)$  exists, and the proof is completed.

\section{ $BF(x, y)$ and $\underline{BF}(x, y)$ with invariant measure}
In this section, we undertake a thorough investigation into the intricate interplay between invariant measures and the functions $BF(x, y)$,  $\underline{BF}(x, y)$. Subsequently, we establish the validity of Theorem 5.4 and Theorem 5.6. The ensuing proposition delineates that in instances where $BF(x, y)=0 $ the measure sets engendered by the dynamical system's trajectories $x$  and  $y$ exhibit identical characteristics.

{\bf Proposition 5.1} Let  $(X, T)$  be a topological dynamical system. If  for any  $x, y \in X $, and $BF(x, y)=0 $, then $M_{B,x}=M_{B,y}$.

{\bf Proof } Given  $x, y \in X $ with  $BF(x, y)=0$. In order to show $ M_{B,x}=M_{B,y}$, we need to prove  $M_{B,x} \subset M_{B,y}$  and  $M_{B,y} \subset M_{B,x} $. In the following, we will prove that  $M_{B,x} \subset M_{B,y}$. Correspondence similary,  $M_{B,y} \subset M_{B,x}$  holds. 

Given  $\mu \in M_{B,x}$, there is a integer subinterval sequence $[m_{r},n_{r}]_{r=1}^{\infty}$  such that for any  $f \in C(X)$, we have
$$
\lim _{r \rightarrow+\infty} \frac{1}{n_{r}-m_{r}} \sum_{k=m_{r}}^{n_{r}-1} f\left(T^{k} x\right)=\int_{X} f \mathrm{~d} \mu .
$$

Given  $f \in C(X)$  and  $\varepsilon>0$. There exists  $\delta=\delta(\varepsilon, f)>0 $  such that whenever  $x_{1}, x_{2} \in X$  with  $d\left(x_{1}, x_{2}\right)<\delta$, we have
$$
\left|f\left(x_{1}\right)-f\left(x_{2}\right)\right|<\varepsilon
$$
Let  $L=\max _{z \in X}\{|f(z)|\}$. Given $ \sigma \in S_{n_{r},m_{r}}$, one has
$$
\begin{aligned}
	\left|\frac{1}{n_{r}-m_{r}} \sum_{k=m_{r}}^{n_{r}-1}  f\left(T^{k} x\right)-\frac{1}{n_{r}-m_{r}} \sum_{k=m_{r}}^{n_{r}-1}  f\left(T^{k} y\right)\right|
	&= \frac{1}{n_{r}-m_{r}}\left|\sum_{k=m_{r}}^{n_{r}-1}\left(f\left(T^{k} x\right)-f\left(T^{\sigma(k)} y\right)\right)\right| \\
	&\leq \frac{1}{n_{r}-m_{r}} \sum_{k=m_{r}}^{n_{r}-1}\left|f\left(T^{k} x\right)-f\left(T^{\sigma(k)} y\right)\right| \\
	&\leq \varepsilon \times \frac{\#\left(\left\{k \in \mathbb{N}^{+} \mid d\left(T^{k} x, T^{\sigma(k)} y\right)<\delta, k \in [m_{r},n_{r}-1]\right\}\right)}{n_{r}-m_{r}}\\
	&+2 L \times \frac{\#\left(\left\{k \in \mathbb{N}^{+} \mid d\left(T^{k} x, T^{\sigma(k)} y\right) \geq \delta,  k \in [m_{r},n_{r}-1]\right\}\right)}{n_{r}-m_{r}} \\
	&\leq \varepsilon+2L \times \frac{\#\left(\left\{k \in \mathbb{N}^{+} \mid d\left(T^{k} x, T^{\sigma(k)} y\right) \geq \delta,  k \in [m_{r},n_{r}-1]\right\}\right)}{n_{r}-m_{r}} .
\end{aligned}
$$
Since
$$
\frac{1}{n_{r}-m_{r}} \sum_{k=m_{r}}^{n_{r}-1}d\left(T^{k} x, T^{\sigma(k)} y\right) \geq \delta \times \frac{\#\left(\left\{k \in \mathbb{N}^{+}\mid d\left(T^{k} x, T^{\sigma(k)} y\right) \geq \delta, k \in [m_{r},n_{r}-1]\right\}\right)}{n_{r}-m_{r}},
$$
we deduce that
$$
\left|\frac{1}{n_{r}-m_{r}} \sum_{k=m_{r}}^{n_{r}-1} f\left(T^{k} x\right)-\frac{1}{n_{r}-m_{r}} \sum_{k=m_{r}}^{n_{r}-1} f\left(T^{k} y\right)\right| \leq \varepsilon+\frac{2 L}{\delta} \times \frac{1}{n_{r}-m_{r}} \sum_{k=m_{r}}^{n_{r}-1} d\left(T^{k} x, T^{\sigma(k)} y\right) .
$$
Thus we get that
$$
\left|\frac{1}{n_{r}-m_{r}} \sum_{k=m_{r}}^{n_{r}-1} f\left(T^{k} x\right)-\frac{1}{n_{r}-m_{r}} \sum_{k=m_{r}}^{n_{r}-1} f\left(T^{k} y\right)\right| \leq \varepsilon+\frac{2 L}{\delta} \times \inf _{\sigma \in S_{n_{r},m_{r}}} \frac{1}{n_{r}-m_{r}} \sum_{k=m_{r}}^{n_{r}-1} d\left(T^{k} x, T^{\sigma(k)} y\right) .
$$
Let  $r \rightarrow+\infty$, then we have
$$
\begin{aligned}
	\limsup _{r \rightarrow+\infty}\left|\frac{1}{n_{r}-m_{r}} \sum_{k=m_{r}}^{n_{r}-1} f\left(T^{k} x\right)-\frac{1}{n_{r}-m_{r}} \sum_{k=m_{r}}^{n_{r}-1} f\left(T^{k} y\right)\right| & \leq \varepsilon+\frac{2 L}{\delta} \times \limsup _{r \rightarrow+\infty} \inf _{\sigma \in S_{n_{r},m_{}}} \frac{1}{n_{r}-m_{r}} \sum_{k=m_{r}}^{n_{r}-1} d\left(T^{k} x, T^{\sigma(k)} y\right) \\
	& =\varepsilon+\frac{2 L}{\delta} \times BF(x, y)=\varepsilon .
\end{aligned}
$$
Let  $\varepsilon \rightarrow 0 $, then we deduce that
$$
\limsup _{r \rightarrow+\infty}\left|\frac{1}{n_{r}-m_{r}} \sum_{k=m_{r}}^{n_{r}-1} f\left(T^{k} x\right)-\frac{1}{n_{r}-m_{r}} \sum_{k=m_{r}}^{n_{r}-1} f\left(T^{k} y\right)\right|=0 .
$$
Thus, we have
$$
\lim _{r \rightarrow+\infty} \frac{1}{n_{r}-m_{r}} \sum_{k=m_{r}}^{n_{r}-1} f\left(T^{k} y\right)=\lim _{r \rightarrow+\infty} \frac{1}{n_{r}-m_{r}} \sum_{k=m_{r}}^{n_{r}-1} f\left(T^{k} x\right)=\int_{X} f \mathrm{~d} \mu,
$$
which implies  $\mu \in M_{B,y}$. Therefore,  $M_{B,x} \subset M_{B,y}$. The proof is completed.

With respect to  $\underline{BF}(x, y) $, we have the following result similar to Proposition 5.1.

{\bf Proposition 5.2} Let  $(X, T)$  be a topological dynamical system. Then  $M_{B,x} \cap M_{B,y} \neq \emptyset$  for any $ x, y \in X$  with  $\underline{BF}(x, y)=0 .$

{\bf Proof }  Given  $x, y \in X$  with  $\underline{BF}(x, y)=0 $. There exists a integer subinterval sequence $[m_{r},n_{r}]_{r=1}^{\infty}$   such that

$$\lim _{r \rightarrow+\infty} \inf _{\sigma \in S_{n_{r},m_{r}}}\frac{1}{n_{r}-m_{r}} \sum_{k=m_{r}}^{n_{r}-1}  d\left(T^{k} x, T^{\sigma(k)} y\right)=0 .$$

Without loss of generality, we can assume that there exists  $\mu \in M_{B,x}$  such that
$$
\lim _{r \rightarrow+\infty} \frac{1}{n_{r}-m_{r}} \sum_{k=m_{r}}^{n_{r}-1}  f\left(T^{k} x\right)=\int_{X} f \mathrm{~d} \mu
$$
holds for any  $f \in C(X) .$
Similar to the process of proving Proposition 5.1, we derive that  $\mu \in M_{B,y} $. Thus,  $M_{B,x} \cap M_{B,y} \neq   \emptyset .$

If  $x \in X$  is a  uniformly generic point of  $(X, T)$, we can strengthen the Proposition 5.1 as follows.

{\bf Proposition 5.3} Let  $(X, T)$  be a topological dynamical system, and  $x, y \in X$. If  $x$  is a uniformly generic point of  $(X, T)$, then  $BF(x, y)=0$  if and only if $ M_{B,x}=M_{B,y} .$

{\bf Proof }  With Proposition 5.1, we only need to show  the sufficiency.

Assume $M_{B,x}=M_{B,y}$. Let $ \mu=\lim _{n-m \rightarrow+\infty} \frac{1}{n-m} \sum_{k=m}^{n-1} \delta_{T^{k} x}$, then  $M_{B,x}=M_{B,y}=\{\mu\} .$  Given  $\varepsilon>0 ,$  and set $ \eta=\frac{\varepsilon}{1+M} $, where  $M=\operatorname{diam}(X) $. According to Lemma 2.3, there are mutually disjoint open sets  $\left\{U_{s}\right\}_{s=1}^{s_{0}} $ such that
$$
\mu\left(\bigcup_{s=1}^{s_{0}} U_{s}\right) \geq 1-\eta \text { and } \operatorname{diam}\left(U_{s}\right) \leq \eta, \forall s=1,2, \cdots, s_{0} .
$$
Combining Lemma 2.2 with Proposition 3.4, we have
$$
\overline{BF}(x, y) \leq \eta \sum_{s=1}^{s_{0}} \mu\left(U_{s}\right)+M\left(1-\sum_{s=1}^{s_{0}} \mu\left(U_{s}\right)\right) \leq \eta+M \eta=\varepsilon .
$$
Let  $\varepsilon \rightarrow 0 $, then we deduce that  $\overline{BF}(x, y)=0 $. This means that $ BF(x, y)=0 .$

Applying Proposition 5.3, we have the following theorem.

{\bf Theorem 5.4} Let  $(X, T)$  be a topological dynamical system. Then  $(X, T)$  is uniquely ergodic if and only if  $N(BF)=X \times X$ if and only if $N(\overline{BF})=X \times X$.

{\bf Proof }  By $N(\overline{BF})=N(BF)$, hence we only need to show $(X, T)$  is uniquely ergodic if and only if  $N(BF)=X \times X$. Assume that  $(X, T)$  is uniquely ergodic and  $\mu $ is the unique ergodic measure. Then for any  $x, y \in X $, we have  $M_{B,x}=M_{B,y}=\{\mu\} $, which implies that  $x$  and  $y$  are  uniformly generic points. According to Proposition 5.3, we derive that  $BF(x, y)=0 $. Thus  $(x, y) \in N(BF)$. Hence  $N(BF)=X \times X .$

Conversely, if  $N(BF)=X \times X $. Let $ \mu_{1}$  and $ \mu_{2}$  be ergodic measures on  $(X, T) $. According to ([4], Theorem 4.28 with page:90), there exist $ x, y \in X $ such that  $M_{B,x}=\left\{\mu_{1}\right\} $ and  $M_{B,y}=\left\{\mu_{2}\right\} $. Since $ N(BF)=X \times X$, we have  $BF(x, y)=0 $. According to Proposition 5.1, we deduce that $ M_{B,x}=M_{B,y}$, which implies $ \mu_{1}=\mu_{2} $. Thus,  $(X, T)$  is uniquely ergodic.

When  $(X, T) $ is a transitive weak Banach mean equicontinuous system, we can deduce that  $N(BF)=   X \times X$. 

Thus according to Theorem 5.4, we have that a transitive weak  Banach mean equicontinuous system is uniquely ergodic as following.

{\bf Corollary 5.5}  Let  $(X, T)$  be a transitive weak Banach mean equicontinuous topological dynamical system. Then  $(X, T) $ is uniquely ergodic.

{\bf Proof }  Let  $x \in X$  be a transitive point of $(X, T)$. Then for any  $y \in X $, there is a subsequence  $\left\{m_{r}\right\}_{r=1}^{+\infty}$  of positive integers  $\mathbb{N}^{+}$ such that  $\lim _{r \rightarrow+\infty} T^{m_{r}} x=y$. Since  $(X, T)$  is weak  Banach mean equicontinuous, we deduce that $ \lim _{r \rightarrow+\infty} \overline{BF}\left(T^{m_{r}} x, y\right)=0$. According to Proposition 5.2, we have that for any $ r \geq 1$,  $\overline{BF}(x, y)=\overline{BF}\left(T^{m_{r}} x, y\right) $. Thus, $ \overline{BF}(x, y)=0 .$

Given  $u, u \in X$. According to Proposition 3.1, we have
$$
\overline{BF}\left(u, v\right) \leq \overline{BF}\left(u, x\right)+\overline{BF}\left(x, v\right)=0
$$
which shows  $\left(u, v\right) \in N(BF)$. Thus $ N(BF)=X \times X $. According to Theorem 5.4, we derive that $ (X, T)$  is uniquely ergodic.

Combining Theorem 5.4 with Proposition 5.2, we can show a new characterization of unique ergodicity based on $ N(\underline{BF}).$

{\bf Theorem 5.6} Let  $(X, T)$  be a topological dynamical system. Then  $(X, T)$  is uniquely ergodic if and only if  $N(\underline{BF})=X \times X$ . 

{\bf Proof }  Assume that  $(X, T)$  is uniquely ergodic. According to Theorem 5,4, we have  $N(BF)=X \times X$. Since  $N(BF) \subset N(\underline{BF})$, we derive that  $N(\underline{BF})=X \times X .$

Conversely, we assume that $ N(\underline{BF})=X \times X .$  Let  $\mu_{1} $ and  $\mu_{2}$  be ergodic measures of  $(X, T)$. According to ([4], Theorem 4.28 with page:90), there are  $x, y \in X$  such that  $M_{B,x}=\left\{\mu_{1}\right\}$  and  $M_{B,y}=\left\{\mu_{2}\right\}$. Since  $N(\underline{BF})=X \times X$, we have  $\underline{BF}(x, y)=0 $. According to Proposition 5.2 , we deduce that $ M_{B,x} \cap M_{B,y} \neq \emptyset $, which implies  $\mu_{1}=\mu_{2}$. Thus,  $(X, T)$  is uniquely ergodic.

Proposition 5.1 shows that  $N(BF)$  is a subset of all point pairs in $ X $ which can generate the same measure set. 

{\bf Theorem 5.7} Let  $(X, T)$  be a topological dynamical system and  $m \in M(X)$. Then the following statements are equivalent:

(1)  $(X, T) $ has uniformly physical measures with respect to  $m$;

(2) $ (m \times m)(N(BF) \cap(Q \times Q))>0 ;$

(3) $ (m \times m)(N(\underline{BF}) \cap(Q \times Q))>0 .$

{\bf Proof }  (2)  $\Leftrightarrow$  (3) According to Theorem 4.2, we have
$$
N(BF) \bigcap(Q \times Q)=N(\underline{BF}) \bigcap(Q \times Q)
$$

(1)  $\Rightarrow $ (2) Let  $\mu$  be a uniformly physical measure of $ (X, T)$  with respect to  $m$. Then $ m(B(\mu))>0$. For any $ x, y \in B(\mu)$, we have  $M_{B,x}=M_{B,y}=\{\mu\}$. Thus  $x, y \in Q $. According to Proposition 5.3 we have  $BF(x, y)=0$, which shows  $(x, y) \in N(BF) \cap(Q \times Q) $. Thus,
$$
B(\mu) \times B(\mu) \subset N(BF) \bigcap(Q \times Q)
$$
Hence we have
$$
(m \times m)(N(BF) \bigcap(Q \times Q)) \geq(m \times m)(B(\mu) \times B(\mu))>0 .
$$
(2)  $\Rightarrow$  (1) There is  $x_{0} \in Q $ such that  $m(N\left(BF, x_{0}\right))>0$. If not, for any  $x \in Q $, we have  $m(N(BF, x))=0 $. Then we derive that
$$
(m \times m)(N(BF) \bigcap(Q \times Q))=\int_{Q} m(N(BF, x)) \mathrm{d} m(x)=0,
$$
which is a contradiction.
Let  $\mu_{x_{0}}$  be the invariant measure generated by  $x_{0}$. Then according to Proposition 5.3, we have  $B\left(\mu_{x_{0}}\right)=   N\left(BF, x_{0}\right) $. Thus
$$
m\left(B\left(\mu_{x_{0}}\right)\right)=m\left(N\left(BF, x_{0}\right)\right)>0,
$$
which implies that $ \mu_{x_{0}}$  is a  uniformly physical measure with respect to  $m .$

\section{Weak Banach mean equicontinuity and Weak  mean equicontinuity}
In this section, we study  $\overline{BF}$-continuity and  $BF$-continuity. Combining the following Proposition 6.4 with Theorem 4.2, we deduce that  $\overline{BF}$-continuity is equivalent to $BF$-continuity. 

{\bf Lemma 6.1}(see [14]) Let  $(X, T)$  be a topological dynamical system. Then  $(X, T)$  is  $\overline{F}$-continuous if and only if  $(X, T)$  is $F$-continuous.

{\bf Lemma 6.2}(see [7]) If a topological dynamical system $(X, T)$ is unique ergodic, then for any $f \in C(X)$ and $x\in X$, 
$
\lim _{n \rightarrow+\infty} \frac{1}{n} \sum_{k=1}^{n} f(T^{k} x)=\lim _{n-m \rightarrow+\infty} \frac{1}{n-m} \sum_{k=m}^{n-1} f(T^{k} x).
$

{\bf Lemma 6.3}(see [14]) Let  $(X, T)$  be an $\overline{F}$-continuous topological dynamical system. Then all the points in $ X$  are  generic points.

{\bf Proposition 6.4} Let  $(X, T)$  be an weak Banach mean equicontinuous topological dynamical system. Then all the points in $ X$  are uniformly generic points.

{\bf Proof }  Given $ x \in X$. For any  $y \in \overline{N(BF, x)}$, there are  $\left\{y_{n}\right\}_{n=1}^{\infty} \subset N(BF, x) $ such that $ \lim _{n \rightarrow \infty} y_{n}=y .$ According to Proposition 3.1, we have
$$
\overline{BF}(x, y) \leq \overline{BF}\left(x, y_{n}\right)+\overline{BF}\left(y_{n}, y\right)=\overline{BF}\left(y_{n}, y\right)
$$
Since  $(X, T) $ is  $\overline{BF}$-continuous, we have  $\lim _{n \rightarrow+\infty} \overline{BF}\left(y_{n}, y\right)=0$. Thus we deduce that
$$
BF(x, y)=\overline{BF}(x, y)=0
$$
which implies $ y \in N(BF, x) $. Hence  $N(BF, x) $ is closed. According to Proposition 3.3, we know  $N(BF, x)$  is an invariant set. Then according to Theorem  5.4, we derive that  $(N(BF, x), T)$  is uniquely ergodic, which implies that all the points in  $N(BF, x) $ are uniformly generic points. In particular,  $x$  is a uniformly generic point.

Since an $BF$-continuous topological dynamical system is $\overline{BF}$-continuous, the following is a direct corollary of Proposition 6.4.  

{\bf Proposition 6.5} Let  $(X, T)$  be an  $BF$-continuous topological dynamical system. Then all the points in $ X$  are uniformly generic points. 

According to Theorem 4.2 and Proposition 6.4, we can get the following theorem.

{\bf Proposition 6.6} Let  $(X, T)$  be a transitive  topological dynamical system. Then $(X, T)$ is weak Banach mean equicontinuous if and only if for every  $x, y \in X$  we have  $\overline{BF}(x, y)=0$ .

{\bf Proof } If  $\overline{BF}(x, y)=0$  for every  $x, y \in X $, then it is clear that  $(X, T)$  is weak Banach mean equicontinuous. Assume that  $(X, T)$  is a transitive  weak Banach mean equicontinuous. Let  $u \in X$  be a point with a dense orbit and take  $y \in X=\overline{\operatorname{Orb}(u, T)} $. And then  we derive that  $\overline{BF}\left(u, y\right)=0$. With the triangle inequality $\overline{BF}(x, y)\leq \overline{BF}(u, x)+\overline{BF}(u, y)$, this implies that  $\overline{BF}(x, y)=0$  for every  $x, y \in X$.
 
{\bf Theorem 6.7} Let  $(X, T)$  be a topological dynamical system. Then  $(X, T)$  is  $\overline{BF}$-continuous if and only if  $(X, T)$  is $BF$-continuous.

{\bf Proof }  We need only to prove that  $\overline{BF}$-continuity implies  $BF$-continuity. Suppose that  $(X, T)$  is  $\overline{BF}$-continuous. According to Proposition 6.4, then we get that all the points in  $X$  are  uniformly generic points. Then according to Theorem 4.2, we have  $BF(x, y)$  exists for any  $x, y \in X$ and $BF(x, y)=\overline{BF}(x, y)$. Thus  $(X, T) $ is  $BF$-continuous.

{\bf Theorem 6.8} Let  $(X, T)$  be a unique ergodic topological dynamical system. Then $(X, T)$  is weak Banach mean equicontinuous if and only if  $(X, T)$  is weak  mean equicontinuous.

{\bf Proof }  It is  clear that a weak Banach mean equicontinuous system  is a  weak  mean equicontinuous system.  Then we need only to indicate that  weak  mean equicontinuity   implies Banach  weak  mean equicontinuity. Suppose that  $(X, T)$  is  weak  mean equicontinuous. According to Lemma 6.3, then $(X, T)$  is unique ergodic and all the points of $ X$  are  generic points. Furthermore with Lemma 6.2, we deduce that all the points of $ X$  are unformly  generic points. Then by Theorem 4.2, for any $x,y\in X$, we have that $BF(x,y)$ exists and $BF(x,y)=F(x,y)$, therefore, $(X, T)$ is $BF$-continuous which implies $(X, T)$  is weak Banach mean equicontinuous. The proof is completed.

{\bf Lemma 6.9}(see [14])  Let  $(X, T)$  be a transitive weak mean equicontinuous topological dynamical system. Then  $(X, T) $ is uniquely ergodic.

{\bf Theorem 6.10} Let  $(X, T)$  be a transitive topological dynamical system. Then $(X, T)$  is weak Banach mean equicontinuous if and only if  $(X, T)$  is weak  mean equicontinuous.

{\bf Proof } By lemma 6.9, a transitive weak mean equicontinuous topological dynamical system is uniquely ergodic, then with the results of Theorem 6.8, we can easily deduce that 
weak  mean equicontinuity  implies Banach  weak  mean equicontinuity. The proof is completed.

{\bf Corollary 6.11} Let  $(X, T)$  be a minimal topological dynamical system. Then $(X, T)$  is weak Banach mean equicontinuous if and only if  $(X, T)$  is weak  mean equicontinuous.

{\bf Proof } We can  easy obtain the result by Theorem 6.10.

{\bf Theorem 6.12} Let  $(X, T)$  be a topological dynamical system. Then  $(X, T)$  is weak Banach mean equicontinuous if and only if the  uniformly time averages  $f_B^{*}$  are continuous for all  $f \in C(X) $.

{\bf Proof }  Assume that  $(X, T) $ is weak Banach mean equicontinuous, then we will prove  $f_B^{*} $ is continuous for any  $f \in C(X) .$

Given  $f \in C(X)$, according to Proposition 6.4 we know $ f_B^{*}(x)$  exists for any $ x \in X$. Fix  $v>0$, then there is $ \varepsilon>0$  such that whenever  $x, y \in X$  with  $d(x, y)<\varepsilon$, we have
$
|f(x)-f(y)|<\frac{v}{2(2 L+1)},
$
where  $L=\max _{x \in X}\{|f(x)|\}$. Given  $x, y \in X$. For any  $\sigma \in S_{n,m} $, we have
$$
\begin{aligned}
	& \left|\frac{1}{n-m} \sum_{k=m}^{n-1} f\left(T^{k} x\right)-\frac{1}{n-m} \sum_{k=m}^{n-1} f\left(T^{\sigma(k)}  y\right)\right| \\
	= & \frac{1}{n-m}\left|\sum_{k=m}^{n-1}\left(f\left(T^{k} x\right)-f\left(T^{\sigma(k)} y\right)\right)\right| \leq \frac{1}{n-m} \sum_{k=m}^{n-1} \mid f\left(T^{k} x-f\left(T^{\sigma(k)} y\right) \mid\right. \\
	\leq & \frac{v}{2(2 L+1)} \times \frac{\#\left(\left\{k \in \mathbb{N}^{+} \mid d\left(T^{k} x, T^{\sigma(k)} y\right)<\varepsilon, k \in [m,n-1]\right\}\right)}{n-m}  \\
	& +2 L \times \frac{\#\left(\left\{k \in \mathbb{N}^{+} \mid d\left(T^{k} x, T^{\sigma(k)} y\right) \geq \varepsilon, k \in [m,n-1]\right\}\right)}{n-m} \\
	\leq & \frac{v}{2(2 L+1)}+2 L \times \frac{\#\left(\left\{k \in \mathbb{N}^{+} \mid d\left(T^{k} x, T^{\sigma(k)} y\right) \geq \varepsilon, k \in[m,n-1]\right\}\right)}{n-m} .
\end{aligned}
$$
Since
$$
\frac{1}{n-m} \sum_{k=m}^{n-1}d\left(T^{k} x, T^{\sigma(k)} y\right) \geq \varepsilon \times \frac{\#\left(\left\{k \in \mathbb{N}^{+} \mid d\left(T^{k} x, T^{\sigma(k)} y\right) \geq \varepsilon, k \in [m,n-1]\right\}\right)}{n},
$$
we deduce that
$$
\left|\frac{1}{n-m} \sum_{k=m}^{n-1} f\left(T^{k} x\right)-\frac{1}{n-m} \sum_{k=m}^{n-1} f\left(T^{k} y\right)\right| \leq \frac{v}{2(2 L+1)}+\frac{2 L}{\varepsilon} \times \frac{1}{n-m} \sum_{k=m}^{n-1} d\left(T^{k} x, T^{\sigma(k)} y\right) .
$$
Thus we have
$$
\left|\frac{1}{n-m} \sum_{k=m}^{n-1} f\left(T^{k} x\right)-\frac{1}{n-m} \sum_{k=m}^{n-1} f\left(T^{k} y\right)\right| \leq \frac{v}{2(2 L+1)}+\frac{2 L}{\varepsilon} \times \inf _{\sigma \in S_{n,m}} \frac{1}{n-m} \sum_{k=m}^{n-1} d\left(T^{k} x, T^{\sigma(k)} y\right) .
$$
Let  $n-m \rightarrow+\infty$, then we have
$$
\limsup _{n-m \rightarrow+\infty}\left|\frac{1}{n-m} \sum_{k=m}^{n-1} f\left(T^{k} x\right)-\frac{1}{n-m} \sum_{k=m}^{n-1} f\left(T^{k} y\right)\right| \leq \frac{v}{2(2 L+1)}+\frac{2 L}{\varepsilon} \times \overline{BF}(x, y) .
$$
This implies
$$
\left|f_B^{*}(x)-f_B^{*}(y)\right| \leq \frac{v}{2(2 L+1)}+\frac{2 L}{\varepsilon} \times \overline{BF}(x, y) .
$$
Since  $(X, T)$  is weak Banach mean equicontinuous, there is  $\delta>0$  such that for any  $x, y \in X $ with  $d(x, y)<\delta$, we have
$$
\overline{BF}(x, y)<\frac{\varepsilon v}{2(2 L+1)}.
$$
Hence, we have that
$$
\left|f_B^{*}(x)-f_B^{*}(y)\right| \leq \frac{v}{2(2 L+1)}+\frac{2 L}{\varepsilon} \times \frac{\varepsilon v}{2(2 L+1)}=\frac{v}{2}<v,
$$
whenever  $x, y \in X$  with  $d(x, y)<\delta$, which implies $ f_B^{*}(x) \in C(X) .$

Conversely, if  $f_B^{*}$  is continuous for any $ f \in C(X) ,$ then we will prove $ (X, T) $ is weak Banach mean equicontinuous.

If  $(X, T) $ is not weak Banach mean equicontinuous, there are  $x \in X,$ $ \varepsilon>0$  and  $\left\{x_{m}\right\}_{m=1}^{\infty} \subset X$  such that  $\lim _{m \rightarrow+\infty} x_{m}=x$  while $\overline{BF}\left(x, x_{m}\right) \geq \varepsilon$. Based on this assumption, we know that $ x $ and $ \left\{x_{m}\right\}_{m=1}^{\infty}$  are uniformly generic points. 

Let  $\mu=\lim _{n-m \rightarrow+\infty} \frac{1}{n-m} \sum_{k=m}^{n-1} \delta_{T^{k} x}$,  $M=\operatorname{diam}(X)$,  $\eta=\frac{\varepsilon}{2(M+5)} .$  Then according to Lemma 2.3, there exist finite mutually disjoint closed subsets  $\left\{\Lambda_{s}\right\}_{s=1}^{s_{0}}$  of  $X$  such that
$$
\mu\left(\bigcup_{s=1}^{s_{0}} \Lambda_{s}\right)>1-\eta \text { and } \operatorname{diam}\left(\Lambda_{s}\right)<\eta, \forall s=1,2, \cdots, s_{0} \text {. }
$$
Take  $\delta=\min _{\delta_{1} \neq s_{2}}\left\{d\left(\Lambda_{s_{1}}, \Lambda_{s_{2}}\right)\right\}, r=\min \left\{\frac{\delta}{5}, \eta\right\} $ and  $\alpha=\frac{\varepsilon}{4 M s_{0}+1} .$ For any  $s \in\left\{1,2, \cdots, s_{0}\right\} $, let $ U_{s}=\left\{y \in X: d\left(y, \Lambda_{s}\right)<r\right\}$  and $ V_{s}=\left\{y \in X: d\left(y, \Lambda_{s}\right)<2 r\right\} $. Then  $\left\{U_{s}\right\}_{s=1}^{s_{0}}$  and  $\left\{V_{s}\right\}_{s=1}^{s_{0}} $ are mutually disjoint open subsets of  $X$  and  $\operatorname{diam}\left(V_{s}\right) \leq 5 \eta$  for any $ s \in\left\{1,2, \cdots, s_{0}\right\}$. We have the following claim:

{\bf Claim 2} For any  $x_{m} \in X$ , there is  $s_{m} \in\left\{1,2, \cdots, s_{0}\right\}$  such that
$$
\liminf _{n -m\rightarrow+\infty} \frac{1}{n-m} \sum_{k=m}^{n-1} \chi_{U_{s_{m}}}\left(T^{k} x\right)>\liminf _{n -m\rightarrow+\infty} \frac{1}{n-m} \sum_{k=m}^{n-1} \chi_{V_{s_{m}}}\left(T^{k} x_{m}\right)+\alpha .
$$
{\bf Proof of Claim 2 } 
 If not, there is  $x_{m} \in X$  such that for any  $s \in\left\{1,2, \cdots, s_{0}\right\}$, we have
$$
\liminf _{n -m\rightarrow+\infty} \frac{1}{n-m} \sum_{k=m}^{n-1} \chi_{U_{s}}\left(T^{k} x\right) \leq \liminf _{n -m\rightarrow+\infty} \frac{1}{n-m} \sum_{k=m}^{n-1} \chi_{V_{s}}\left(T^{k} x_{m}\right)+\alpha .
$$
Given  $s \in\left\{1,2, \cdots, s_{0}\right\} $. According to Lemma 2.1, we have that
$$
\liminf _{n -m\rightarrow+\infty} \frac{1}{n-m} \sum_{k=m}^{n-1}\chi_{U_{s}}\left(T^{k} x\right) \geq \mu\left(U_{s}\right).
$$
Then, we deduce that
$$
\begin{aligned}
	\liminf _{n -m\rightarrow+\infty} \frac{1}{n-m} \sum_{k=m}^{n-1} \chi_{V_{s}}\left(T^{k} x_{m}\right) & \geq \liminf _{n -m\rightarrow+\infty} \frac{1}{n-m} \sum_{k=m}^{n-1} \chi_{U_{s}}\left(T^{k} x\right)-\alpha \\
	& \geq \mu\left(U_{s}\right)-\alpha .
\end{aligned}
$$
As $ U_{s} \subset V_{s}$, we have that
$$
\liminf _{n -m\rightarrow+\infty} \frac{1}{n-m} \sum_{k=m}^{n-1} \chi_{V_{s}}\left(T^{k} x\right) \geq \liminf _{n -m\rightarrow+\infty} \frac{1}{n-m} \sum_{k=m}^{n-1} \chi_{U_{s}}\left(T^{k} x\right) \geq \mu\left(U_{s}\right) .
$$
Then According to Proposition 3.4, we derive that
$$
\begin{aligned}
	\overline{BF}\left(x, x_{m}\right) & \leq 5 \eta \sum_{s=1}^{s_{0}}\left(\mu\left(U_{s}\right)-\alpha\right)+M\left(1-\sum_{s=1}^{s_{0}}\left(\mu\left(U_{s}\right)-\alpha\right)\right) \\
	& \leq 5 \eta+M \eta+M s_{0} \alpha \leq \frac{3 \varepsilon}{4}.
\end{aligned}
$$
which is a contradiction. The proof of Claim is completed.

According to Claim, there is  $s_{m_{0}} \in\left\{1,2, \cdots, s_{0}\right\}$  and a subsequence  $\left\{x_{m_{p}}\right\}_{p=1}^{\infty}$  of  $\left\{x_{m}\right\}_{m=1}^{\infty}$  such that for any  $p \in \mathbb{N}^{+} $, we have
$$
\liminf _{n -m\rightarrow+\infty} \frac{1}{n-m} \sum_{k=m}^{n-1} \chi_{U_{s_{m_{0}}}}\left(T^{k} x\right)>\liminf _{n -m\rightarrow+\infty} \frac{1}{n-m} \sum_{k=m}^{n-1} \chi_{V_{s_{m_{0}}}}\left(T^{k} x_{m_{p}}\right)+\alpha  \quad(6. 5)
$$
Take $ f \in C(X)$  such that $ 0 \leq f \leq 1$  and
$$
\left.f\right|_{\overline{U}_{s_{m_{0}}}}=1,\left.\quad f\right|_{V_{s_{m_{0}}}^{c}}=0$$
Then we have that
$$
\liminf _{n -m\rightarrow+\infty} \frac{1}{n-m} \sum_{k=m}^{n-1} \chi_{U_{s_{m_{0}}}}\left(T^{k} x\right) \leq \liminf _{n -m\rightarrow+\infty} \frac{1}{n-m} \sum_{k=m}^{n-1} f\left(T^{k} x\right)=f_B^{*}(x)
$$
and
$$
\liminf _{n -m\rightarrow+\infty} \frac{1}{n-m} \sum_{k=m}^{n-1} \chi V_{s_{m_{0}}}\left(T^{k} x_{m_{p}}\right) \geq \liminf _{n -m\rightarrow+\infty} \frac{1}{n-m} \sum_{k=m}^{n-1} f\left(T^{k} x_{m_{p}}\right)=f_B^{*}\left(x_{m_{p}}\right) .
$$
Thus, we deduce that
$$
f_B^{*}(x) \geq f_B^{*}\left(x_{m_{p}}\right)+\alpha,
$$
which implies  $f_B^{*}(x) \notin C(X)$. This is a contradiction. 

Therefore, $ (X, T)$  is weak Banach mean equicontinuous.

{\bf Lemma 6.13}(see [14]) Let  $(X, T)$  be a topological dynamical system. Then  $(X, T)$  is weak  mean equicontinuous if and only if the  time averages  $f^{*}$  are continuous for all  $f \in C(X) $.

By Theorem 6.8, Theorem 6.11, Theorem 6.12 and Lemma 6.13, we can obtain the results as follows.

{\bf Corollary 6.14} Let  $(X, T)$  be a unique topological dynamical system. Then  the  uniformly time averages  $f_B^{*}(x)$  are continuous for all  $f \in C(X) $ if and only if the  time averages  $f^{*}$  are continuous for all  $f \in C(X)$.

{\bf Corollary 6.15} Let  $(X, T)$  be a transitive topological dynamical system. Then  the  uniformly time averages  $f_B^{*}(x)$  are continuous for all  $f \in C(X) $ if and only if the  time averages  $f^{*}$  are continuous for all  $f \in C(X)$.

{\bf Corollary 6.16} Let  $(X, T)$  be a minimal topological dynamical system. Then  the  uniformly time averages  $f_B^{*}(x)$  are continuous for all  $f \in C(X) $ if and only if the  time averages  $f^{*}$  are continuous for all  $f \in C(X)$.

\bigskip \noindent{\bf Acknowledgement}.
The research was supported by NSF of China (No. 11671057) and NSF of Chongqing (Grant No. cstc2020jcyj-msxmX0694).


 \bibliographystyle{amsplain}

\end{document}